\documentclass[12pt]{article}

\usepackage{amsmath}
\usepackage{amsfonts}
\usepackage{amssymb}

\usepackage{mathstyle}
\usepackage{flexisym}
\usepackage{breqn}

\usepackage{srcltx}

\usepackage{pgf}
\input xy
\xyoption{all}

\usepackage{hyperref}

\newtheorem{theo}{Theorem }[section]
\newtheorem{lemm}{Lemma}[section]

\newtheorem{coro}{Corollary}[section]

\newtheorem{defi}{Definition}[section]

\newtheorem{defii}{Definition}

\newtheorem{rema}{Remark}[section]

\newtheorem{exa}{Example}[section]

\def\corank{\mathop{\mathrm{corank}}}

\def\codim{\mathrm{codim}}

\def\Id{\mathrm{Id}}

\def\inv{\mathrm{inv}}

\def\g{\mathfrak{g}}
\def\gl{\mathfrak{gl}}

\def\Lie{\mathrm{Lie}}

\def\rank{\mathop{\mathrm{rank}}}
\def\s{\mathfrak{s}}
\def\sl{\mathfrak{sl}}
\def\so{\mathfrak{so}}

\def\Sing{\mathop{\mathrm{Sing}}}

\def\Span{\mathop{\mathrm{Span}}}

\def\Trace{\mathop{\mathrm{Trace}}}

\def\Vect{\mathop{\mathrm{Vect}}}

\def\K{\mathbb{K}}

\def\R{\mathbb{R}}

\def\L{\mathcal{L}}
 \def\C{\mathbb{C}}

\def\al{\alpha}

\def\d{\partial}

\def\la{\lambda}

\def\om{\omega}

\def\qed{$\square$}

\title{On linear-quadratic Poisson pencils on trivial central extensions of semisimple Lie algebras}
\author{Andriy Panasyuk\\
and\\
Vsevolod Shevchishin\\
\\
Faculty of Mathematics and Computer Science\\
University of Warmia and Mazury\\
ul. S\l oneczna 54, 10-710 Olsztyn, Poland\\
\\
panas@matman.uwm.edu.pl\\
shevchishin@googlemail.com
}
\date{}

\begin{document}
\bibliographystyle{plain}

\maketitle

\begin{abstract}
The paper is devoted to quadratic Poisson structures compatible with the canonical linear Poisson structures on  trivial 1-dimensional central extensions of semisimple Lie algebras. In particular, we develop the general theory of such structures and study related families of functions in involution. We also show that there exists a 10-parametric family of quadratic Poisson structures on $\gl(3)^*$ compatible with the canonical linear Poisson structure and containing the 3-parametric family of quadratic bivectors recently introduced by Vladimir Sokolov. The involutive family of polynomial functions related to the corresponding Poisson pencils contains the hamiltonian
 of the polynomial form of the elliptic Calogero--Moser system.
\end{abstract}
\section*{Introduction}

It is well-known that  the  bihamiltonian property proved to be very useful in the study of integrable systems. Since its invention by Magri  \cite{m} in the context of the KdV equation a second hamiltonian structure was discovered for the majority of important examples. By a bihamiltonian property we mean the existence of second Poisson structure, compatible with the first one (canonically defined), such that the vector field defining the evolution of the system is hamiltonian in two ways, i.e. with respect to both the Poisson structures.

Recall that Poisson bivectors\footnote{We will use the term ``polyvector'' instead of ``polyvector field'' for short, with the exception of vector fields.} $\pi_1$ and $\pi_2$ on a smooth manifold $M$ are called \emph{compatible } if $\pi_1+\pi_2$ is a Poisson bivector, or, equivalently, $[\pi_1,\pi_2]=0$, where $[,]$ stands for the Schouten bracket. Any pair $(\pi_1,\pi_2)$ of compatible Poisson bivectors generates a two-dimensional linear space of Poisson bivectors $\{\la_1\pi_1+\la_2\pi_2\}_{(\la_1,\la_2)\in\K^2}$ which is called a \emph{Poisson pencil} (here $\K$ is the ground field).

By linear, quadratic etc. polyvectors  on a vector space we mean polyvectors with {\em homogeneous }linear, quadratic etc. coefficients with respect to some linear system of coordinates.

We will say that a Poisson pencil is \emph{linear-quadratic} if its generators are linear and quadratic bivectors on some vector space respectively. The cases of ``linear-constant'' and ``linear-linear'' Poisson pencils are already very important, see for instance \cite{bols'}, \cite{pBi-Lie} and references therein. The linear-quadratic Poisson pencils are the main objects of this paper.

There are basically two classes of linear-quadratic Poisson pencils known in the literature.
The first one can be described as follows (\cite{doninGurevich}, \cite{gurevichPanyushev}, \cite{gurevichRubtsov}). Let $r\in\g\wedge\g$ be a skew-symmetric solution to the classical Yang--Baxter equation on a Lie algebra $\g$, i.e. an "algebraic Poisson bivector"
on $\g$. Assume that a linear representation of $\g$ in a vector space $V$ is given and $\xi_1,\ldots,\xi_n$ are the fundamental vector fields of this representation corresponding to a basis $e_1,\ldots,e_n$ of $\g$. Then, if $r=r^{ij}e_i\wedge e_j$, the bivector $\pi_2:=r^{ij}\xi_i\wedge \xi_j$ is a quadratic Poisson bivector  on $V$. If, moreover, another bivector $\pi_1$ on $V$ is given, invariant under the action of $\g$, the bivectors $\pi_1$ and $\pi_2$ are compatible. In particular, if $V=\g^*$ with the canonical Lie--Poisson bivector $\pi_1$, and the representation is the coadjoint one, we get a linear-quadratic Poisson pair $(\pi_1,\pi_2)$, where $\pi_2=r^{ij}\pi_1(x_i)\wedge
\pi_1(x_j)$ (here $x_i$ are the corresponding coordinates on $\g^*$).

The second class concerns a specific situation when a quadratic Poisson bivector $\pi_2$ is given on a vector space $V$ and there exists a point $a\in V$, $a\not=0$, such that $\pi_2(a)=0$ (not for every quadratic Poisson bivector such points exist). Then the \emph{linearization} of $\pi_2$ at $a$ is a linear Poisson bivector which is automatically compatible with $\pi_2$.

One can find numerous examples of such a situation in the literature. One class of examples is related to Poisson--Lie groups. Given such a group $(G,\pi)$, its Poisson tensor $\pi$ necessarily vanishes at the neutral element $e\in G$ and there is a correctly defined linearization $\pi_1$ of $\pi$ defined on the linear space $T_eG$ (it coincides with the $r$-matrix bracket on $\g=\Lie(G)$ of the corresponding tangent Lie bialgebra). Sometimes $G$ can be embedded in a vector space which can be identified with $T_eG$ (as in the case of the group $GL(n,\K)$) and the two bivectors appear to be compatible (\cite{sts1},\cite[\S 1--2]{karasevMaslov}).

Another important class of examples are the so-called elliptic quadratic Poisson bivectors $\pi_2$ on $\gl(n,\C)$ which are the quasiclassical limits of the antisymmetrization of the  Feigin--Odesskii algebras (which coincide with the Sklyanin algebras for $n=2$). Here $\pi_2(I)=0$, where $I$ is the unity matrix, and the linearization of $\pi_2$ at $I$ is the standard Lie--Poisson bivector $\pi_1$ on $\gl(n,\C)\cong\gl^*(n,\C)$. The resulting linear-quadratic Poisson pencil is related to the so-called elliptic rotator bihamiltonian integrable system (\cite{khesinLO}). The coefficients of the tensors $\pi_2$ are written by means of elliptic constants.

Recently,   V. Sokolov introduced a family of quadratic Poisson bivectors $\pi_2$ on the space $\gl(3,\C)^*$  which
have the same property: $\pi_2(I)=0$ and the linearization of $\pi_2$ at $I$ coincides with the standard Lie--Poisson bivector $\pi_1$ on $\gl(3,\C)\cong\gl^*(3,\C)$  \cite{SokolovEllCM}, \cite{SokolovObzor}. This family is also related to an elliptic curve but on the contrary to the above case the coefficients of $\pi_2$ are written explicitly (with no use of elliptic constants). They depend on the invariants $g_2,g_3$ of the elliptic curve, however one can easily check that  these numbers can take any values and in fact the Sokolov family determines a 3-parametric family of  quadratic Poisson bivectors (the third component corresponds to the free term).

The involutive family of polynomial functions canonically associated with the Poisson pair $(\pi_1,\pi_2)$ (see Definition \ref{canoo})  contains the hamiltonian of the classical elliptic Calogero--Moser system. More precisely, in \cite{MatushkoSok} the hamiltonian operator of the so-called quantum elliptic Calogero--Moser system for 3 particles is brought to a polynomial form. Moreover, a system of commuting elements of the universal enveloping algebra $U(\gl(3))$ is built such that under a specific representation these elements pass to the differential operators:  the second order quantum hamiltonian and its quantum integrals of third order. The above mentioned involutive family of  functions is constructed by considering the symbols of the commuting elements of $U(\gl(3))$.

Notice that the Sokolov quadratic bivectors, as well as the  elliptic quadratic Poisson bivectors \cite{khesinLO}, are \emph{centrally linearizable}  in the sense of the following definition.
\begin{defii}\label{efi}\rm Let $\s$ be a semisimple Lie algebra and  let $\g=\s\oplus\K$ be its one-dimensional trivial central extension. A quadratic Poisson bivector
$\pi_2$ on $\g^*$  is {\em centrally linearizable} (CL for short) if
\begin{enumerate}\item $\pi_2(a)=0$ for some $a\in \K^*\subset\g^*=\s^*\oplus\K^*$, $a\not=0$;
\item the linearization of $\pi_2$ at $a$ is equal to the canonical Lie--Poisson bivector $\pi_1$ on $\g^*$.
\end{enumerate}
\end{defii}

It is well-known that condition 2 implies that the bivectors $\pi_1$ and $\pi_2$ are compatible; we also call \emph{centrally linearizable} the corresponding Poisson pencil.

The interest of the authors to this topic started from the observation that the Sokolov 3-parametric family of CL quadratic Poisson bivectors can be included into a 10-parametric family of CL quadratic  Poisson bivectors. The aim of this paper is to develop a theory of CL Poisson structures  and to prove the above  fact along with the explicit calculation of the canonically defined family of integrals and the discussion of possibility of generalization of these results to other Lie algebras.

We show in Section \ref{adm} that  CL Poisson bivectors are in a 1--1 corrspondence with  pairs $(X,q)$, where $X$ is a quadratic vector field and $q$ a quadratic polynomial on $\s^*$, satisfying a couple of nonlinear identities (\ref{eqbasic}), (\ref{eqbasic2}), which we call the \emph{basic identities} (we prove that in fact (\ref{eqbasic}) implies (\ref{eqbasic2}) for higher dimensional Lie algebras, Theorem \ref{th2}). The following observation was crucial: all three quadratic vector fields $X$ corresponding to the components of the Sokolov family modulo inessential quadratic hamiltonian vector fields  belong to a 10-dimensional $\sl(3,\C)$-irreducible invariant subspace $W^*$ in the space of all quadratic vector fields on $\s^*=\sl(3,\C)^*$ (cf. Remark \ref{sokk}).  We next are able, given any quadratic vector field $X\in W^*$, to find a quadratic polynomial $q$ such that the pair $(X,q)$ satisfies the basic identities.

We call a CL quadratic Poisson bivector $\pi$ on $\g^*=(\s\oplus\K)^* ${\em monogenic} if the corresponding quadratic vector field $X$ belongs to an $\s$-irreducible invariant subspace $\mathcal{W}$ in the space of all the quadratic vector fields on $\s^*$ such that for any $Y\in \mathcal{W}$ there exists a quadratic polynomial $q$ satisfying the basic identities along with $Y$ (i.e. the CL bivector $\pi$ can be included into a $\dim \mathcal{W}$-parametric family of CL Poisson bivectors). In particular, the CL Poisson bivectors from the 10-dimensional family mentioned are monogenic.

A natural question arise: are there other monogenic CL Poisson bivectors (on $\gl(3)$ or other Lie algebras)? We propose the following approach for looking for monogenic CL bivectors. Try, given a zero weight quadratic vector field $X$ belonging to some irreducible $\s$-module, to find a quadratic polynomial (necesssarily being of zero weight) satisfying the basic identities along with $X$. We conjecture that  the corresponding CL bivector will be monogenic and the rest monogenic bivectors corresponding to  this module could be found by equivariance, in analogy to the case of the module $W^*$. Note that the number of parameters of the zero weight spaces is essentially less than that of the whole module, which simplifies the calculations. On the other hand, one cannot use for this purpose other weight spaces, since, for instance, the quadratic polynomial $q$ being the solution of the basic identities for the highest weight vector $X$ in $W^*$ is trivial.

We applied this approach to all the irreducible submodules of the 288-dimensional module $S^2\s\otimes\s^*$  of quadratic vector fields on $\s^*=\sl(3,\C)^*$ (cf. Section \ref{s10}) and found that apart the 10-dimensional irreducible module $W^*$ and its dual $W=(W^*)^*$  (cf. footnote \ref{secc}) there is only one more module leading to monogenic bivectors. This module is isomorphic to $\s$ itself and consists of the vector fields of the form $xE$, where $x$ is a linear function on $\s^*$ and $E$ is the Euler vector field on $\s^*$. The corresponding quadratic polynomial $q$ is equal to $x^2/2$ and in fact this solution to the basic identities, which we call {\em tautological}, is a universal one, i. e. the construction works for any Lie algebra $\s$ (even not semisimple one), see Theorem \ref{tauto}.

We conjecture that  it is possible  to find (nontautological) monogenic  CL Poisson bivectors for trivial 1-dimensional central extensions of other semisimple Lie algebras.  In general one may expect the existence of nonmonogenic CL quadratic bivectors, i.e. such that the corresponding vector field $X$  consists of several components belonging to different irreducible modules (it seems
 that the  elliptic quadratic Poisson bivector on $\gl(n,\C)$ \cite{khesinLO} is of this type).
Thus the question of classification of all the CL Poisson bivectors looks very complicated even for low-dimensional Lie algebras. Another important question is related to  quantization of the CL quadratic Poisson bivectors (for instance in the deformation quantization sense) and of the corresponding involutive families of functions. We conjecture that the functions in involution corresponding to our 10-dimensional family of linear-quadratic Poisson pencils on $\gl^*(3)$ are the symbols of  commuting elements in $U(\gl(3))$ in complete analogy with the Sokolov case.

Finally, we would like to underline the exceptional role of the CL Poisson quadratic bivectors among other quadratic Poisson bivectors compatible with the linear ones. Indeed, any quadratic Poisson bivector $\pi_2$ compatible with the canonical linear Poisson bivector $\pi_1$  on $\g^*=\s^*$ ($\s$ semisimple) or $\g^*=(\s\oplus\K)^*$ is of the form $[X,\pi_1]$, where $X$ is a quadratic vector field satisfying $[[X,\pi_1],[X,\pi_1]]=0$ (see Section \ref{cobo}). It is not hard to find examples of such vector fields $X$ on $\g^*$, for instance so are the highest vectors of some of irreducible submodules in $S^2\s\otimes\s^*$. However, there is no guarantee that the canonical families of commuting functions related to the Poisson pair $(\pi_1,\pi_2)$ contain quadratic functions  (which could serve as ``physically reasonable'' hamiltonians). On the other hand, in the case of CL Poisson pencils  the quadratic hamiltonians are in a sense incorporated in the theory in the form of the polynomials $q$, which as Theorem \ref{th4a} shows are included in the canonical family of commuting functions.

Let us review the content of the paper in more detail. In Section \ref{PP} we recall some very basic facts about the related with Poisson pencils canonical families of functions in involution, in particular about the Magri--Lenard chains. Section \ref{cobo} is devoted to the proof of the fact that any  quadratic Poisson bivector $\pi_2$ compatible with the canonical linear Poisson bivector $\pi_1$  on $\s^*$ ($\s$ semisimple) or $\g^*=(\s\oplus\K)^*$ is of the form $[X,\pi_1]$, where $X$ is a quadratic vector field. We also prove that, given a linear-quadratic Poisson pencil, there exists a Magri--Lenard chain starting from any Casimir function of the linear Poisson bivector. In Section \ref{adm} we develop the general theory of CL quadratic Poisson pencils, in particular we discuss the matters of Casimir functions, the canonical involutive families of functions and their completeness. Elementary examples of CL pencils, including the Sklyanin bracket, the above mentioned tautological solution of the basic identities and a generalization of brackets related to a constant solution of the Yang--Baxter equation, appear in Section \ref{taut}. Section \ref{s10} introduces two irreducible modules: a 10-dimensional one in the space of quadratic vector fields (the above mentioned module $W^*$) and a 27-dimensional submodule $U$ in the space of quadratic polynomials on $\sl(3)^*$ (the polynomials $q$ solving the basic identities along with $X\in W^*$ belong to $U$). In Section \ref{ser} we show the existence of a 10-parametric family of CL quadratic Poisson bivectors on $\gl(3)^*$ and we discuss the corresponding Casimir functions and the canonical families of functions in involution in Section \ref{cas1}. The explicit formulas for the corresponding vector fields and functions depending on 10 parameters are given in Appendix.

In this paper we permanently use the standard properties of the Schouten bracket \cite{marle2}. A part of the  calculations was made using the Maple software as well as the LiE program (http:wwwmathlabo.univ-poitiers.fr/$\sim$maavl/LiE/). We tried to find theoretical explanations or proofs for all the results, however we failed in some cases, in which the proof is done by  direct computer check.


\section{Preliminaries on Poisson pencils}
\label{PP}

If $M$ is a real or complex analytic manifold, ${\cal E}(M)$ will stand
for the space of  analytic  functions on $M$ in the corresponding category.
We say that a Poisson pencil  $\{\la_1\pi_1+\la_2\pi_2\}$ on a  manifold $M$ is of \emph{Gelfand--Zakharevich type} if $\rank(\la_1\pi_1+\la_2\pi_2)|_x<\dim M$ for any $(\la_1,\la_2)$ and any $x\in M$.  In this paper we consider only Poisson pencils of Gelfand--Zakharevich type with the following additional property: there exist an open dense set $U \subset M$ and a collection of functions $C_1,\ldots,C_r\in{\cal E}(M)$, where $r=\min_{x\in M}\corank(\la_1\pi_1+\la_2\pi_2)|_x$, which are functionally independent on $U$ and are global Casimir functions of the bivector $\pi_1$, i.e. $\pi_1(C_i):=[\pi_1,C_i]=0$ on $M$ (cf. footnote \ref{fn}). Such a property holds for instance for the Lie--Poisson structure $\pi_1$ on $\g^*$, where $\g$ is a semisimple or reductive Lie algebra.

\begin{defi}\label{defML} \rm Let $\{\la_1\pi_1+\la_2\pi_2\}$ be a Poisson pencil.
An infinite sequence of functions $f_0,f_1,\ldots\in{\cal E}(M)$, is called a \emph{Magri--Lenard chain} if it satisfies the relations
\begin{equation}\label{eq04}
 \pi_1(f_0)=0,\pi_1(f_1)+\pi_2(f_0)=0,\ldots,\pi_1(f_{i+1})+\pi_2(f_i)=0,\ldots.
\end{equation}
\end{defi}
Note that  $f_0,f_1,\ldots$ is a Magri--Lenard chain if and only if  the formal series $f(\la)=f_0+\la f_1+\cdots$ is a formal Casimir function of the bivector $\pi_1+\la\pi_2$ (cf. the notions of ``anchored formal $\la$-family'' and ``anchored Lenard scheme'' in \cite[Sect. 10]{gz4}). Obviously, if $f_0,f_1,\ldots$ is a Magri--Lenard chain such that $f_k\not=0$ and $f_i=0$ for $i>k$, then $f(\la)$ is a Casimir function of $\pi_1+\la\pi_2$ and $f_k$ is a Casimir function of the bivector $\pi_2$.

\begin{rema}\label{invv}\rm By \cite{gz4}[Prop. 10.7], given two Magri--Lenard chains $f_0,f_1,\ldots$ and $g_0,g_1,\ldots$, all their members commute with each other with respect to the Poisson bracket corresponding to any $\la_1\pi_1+\la_2\pi_2$.
\end{rema}

\begin{defi}\label{canoo} \rm
The family of functions functionally generated by all the members of all the Magri--Lenard chains starting from all the Casimir functions of the bivector $\pi_1$ for which such chains exist will be called the \emph{canonical involutive family of functions} related to the pencil $\{\la_1\pi_1+\la_2\pi_2\}$.
\end{defi}

\begin{rema}\label{remmm}\rm
One can consider also {\em finite} sequences of functions $f_0,\ldots,f_{i+1}$ satisfying relations (\ref{eq04}) with the property that there does not exist a function $f_{i+2}$ satisfying the relation $\pi_1(f_{i+2})+\pi_2(f_{i+1})=0$. They correspond to the Jordan blocks in the so-called Jordan--Kronecker decomposition of the pair of skew-symmetric forms $\pi_1|_x,\pi_2|_x:T_x^*M\to T_xM$, $x \in M$
\cite{gz4}, \cite{pSymmetries}, \cite{bolsSing}. The members of such families fail to commute in general.

On the contrary, the Magri--Lenard chains correspond to the Kronecker blocks in the Jordan--Kronecker decomposition (the Gelfand--Zakharevich type structures necessarily contain at least one such block). This shows that the family of functions generated by the members of the Magri--Lenard chains is canonical indeed, i.e. does not depend on the choice of basis $\pi_1,\pi_2$ in the Poisson pencil.

Finally, it is worth mentioning that there is a class of Poisson pencils called \emph{Kronecker} such that the Jordan--Kronecker decomposition of the pair $\pi_1|_x,\pi_2|_x$ contains only Kronecker blocks for any $x$ from some open and dense set in $M$. In the case of Kronecker Poisson pencils the canonical family of functions is  \emph{complete} \cite[Prop. 11.7, 11.11]{gz4}, i.e. cuts a lagrangian foliation on generic symplectic leaves of $\la_1\pi_1+\la_2\pi_2$ for nontrivial $(\la_1,\la_2)$.
\end{rema}

\section{Coboundary Poisson bivectors and linear-quadratic Poisson pairs}
\label{cobo}

We say that a Poisson bivector $\eta$ is of {\em coboundary type} if there exist a Poisson bivector $\pi$ and a vector field $X$ such that $\eta=[X,\pi]=\L_X\pi$ (here $[,]$ stands for the Schouten bracket and $\L_X$ for the Lie derivative). In particular, $\eta$ is a coboundary in the sense of the Poisson cohomology of $\pi$ \cite[Sect. 3.6]{wcs}.

Note that $\eta$ and $\pi$ are automatically compatible, i.e. $[\eta,\pi]=0$, since $[\eta,\eta]=0$ and due to the graded anticommutativity and Jacobi identity for the Schouten bracket. This observation allows to construct Poisson pencils by means of a pair $(X,\pi)$ consisting of a vector field and a Poisson bivector provided condition $[[X,\pi],[X,\pi]]=0$ is satisfied.

Now let $\pi_1$ be a {\em linear} Poisson bivector on $\g^*$, where $\g$ is the corresponding Lie algebra, and $X$ be a {\em quadratic} vector field on $\g^*$. Then $[X,\pi_1]$ is a quadratic bivector.

\begin{lemm}\label{lem}
Let $\g$ be a Lie algebra. Assume that $H^i(\g,S^2(\g))=0$, $i=1,2$,  and let $\pi_1$ be the corresponding Lie--Poisson bivector. Then
\begin{enumerate}\item
 each quadratic Poisson bivector on $\g^*$ compatible with $\pi_1$ is of the form $[X,\pi_1]$ for some quadratic vector field $X$;
\item the vector field $X$ is defined up to the addition of a Poisson cocycle (i.e. a vector field $Z$ such that $[Z,\pi_1]=0$), which necessarily is a hamiltonian vector field $\pi_1(f)$\footnote{\label{fn} Note that $\pi_1(q)=[\pi_1,q]$ is the contraction of $\pi_1$ and $dq$ over the second index, cf. \cite[Form. (20)]{marle2}, so our notion differs by sign from the common notion of the hamiltonian vector field (the contraction over the first index).} for some quadratic polynomial $f$.
\end{enumerate}
\end{lemm}

\noindent The proof follows from the fact that the Poisson cohomology with respect to the Lie--Poisson structure $\pi_1$ restricted to the quadratic polyvectors is naturally identified with the Chevalley--Eilenberg cohomology $H^\bullet(\g,S^2\g)$ of the Lie algebra $\g$ with coefficients  in the module $S^2\g$ \cite[Lemma 4.5]{dufourZung}. \qed

\begin{rema}\label{rema}\rm By the Whitehead lemmas $H^i(\g,M)=0$, $i=1,2$, for any semisimple Lie algebra and any finite-dimensional module, in particular the lemma holds for a semisimple $\g$. For the one-dimensional trivial central extension $\g=\s\oplus \K $ of a semisimple Lie algebra $\s$ over the field $\K$ the situation is as follows : $H^1(\g,S^2(\g))=\K^2$ and $H^2(\g,S^2(\g))=0$. In particular, Item 1 of Lemma \ref{lem} remains true and Item 2 has to be corrected, as we explain below.

By \cite{FeiginFuchs}[Ch. 3,\S 1] we have $H^i(\K,\K)=\K$ for $i=0,1$ and $H^i(\K,\K)=0$ for $i>1$. Hence $H^1(\g,\K)=H^1(\K,\K)\otimes H^0(\s,\K)=\K$, $H^2(\g,\K)=\sum_{i+j=2}H^i(\K,\K)\otimes H^j(\s,\K)=0$ and $H^1(\g,S^2(\g))=H^1(\g,\K)\otimes \inv S^2(\g)=\K^2$, $H^2(\g,S^2(\g))=H^2(\g,\K)\otimes \inv S^2(\g)=0$, where $\inv S^2(\g)$ stands for the space of invariant quadratic polynomials on $\g^*$ which is 2-dimensional (generated by $x_0^2$ and $B$, where  $x_0$ is the coordinate corresponding to the centre and $B$ is induced by the Killing form  of $\s$)). In particular,  the vector field $X$ is defined up to the addition of a Poisson cocycle, which is the sum of a hamiltonian vector field $\pi_1(f)$ and the vector field of the form $(ax_0^2 + bB)\frac{\d }{\d x_0}$, $a,b\in\K$. Analogously one can prove that any linear vector field being a $\pi_1$-cocycle is a sum of a hamiltonian vector field $\pi_1(g)$ with linear $g$ and the vector field proportional to $x_0\frac{\d }{\d x_0}$.
\end{rema}

\begin{lemm}\label{lem22}
Let $\g$ be a Lie algebra such that $H^2(\g,S^r(\g))=0$, $r=3,4,\ldots$ (in particular $\g$ can be semisimple or the one-dimensional trivial central extention of a semisimple Lie algebra). Let $\pi_1$ be the corresponding Lie--Poisson bivector and let $\pi_2=[X,\pi_1]$ be a quadratic Poisson bivector compatible with $\pi_1$. Then
\begin{enumerate}\item
 there exists a formal trivializing map $\Phi(t):\g\to\g$, $\Phi(t)=\Id_\g+tZ_1+t^2Z_2+\cdots$, $\Phi_*(\pi_1+t\pi_2)=\pi_1$, where $Z_j:\g\to\g$ is a homogeneous polynomial map of degree $j+1$, i.e. for any $j$ the following equality holds
 $$
 (\Id_\g+tZ_1+\cdots +t^{j}Z_i)_*(\pi_1+t\pi_2)=\pi_1+o(t^j);
 $$
  here $o(t^j)$ stands for the terms of degree $>j$ in $t$ beinng bivectors of degree $>j+1$;
  \item if $\{x_i\}$ is a system of linear coordinates on $\g$ and $X=X_i(x) \frac{\d }{\d x_i}$, then $Z_1$ is given by $x_i\mapsto X_i(x)$;
\item if $f$ is a Casimir function of the bivector $\pi_1$, $\pi_1(f)=0$, then $\Phi(t)^*f=f+tf_1+t^2f_2+\cdots$ is a formal Casimir function of the bivector $\pi_1+t\pi_2$, i.e. the functions $f_0:=f,f_1,\ldots$ form a Magri--Lenard chain (see Definition \ref{defML}),
    where in particular,
    \begin{equation}\label{ca}
    f_1=Xf.
    \end{equation}
   \end{enumerate}
\end{lemm}

\noindent Let  $Y$ be any vector field,  $\eta$  any bivector, and let $\{,\}$ denote  the skew-symmetric bracket on functions corresponding to the bivector $\eta$, $\{f,g\}=[[\eta,g],f]$.
We exploit the formula for the bracket on smooth functions corresponding to the bivector $[Y,\eta]$:
$$
\{f,g\}_{[Y,\eta]}=Y\{f,g\}-\{Yf,g\}-\{f,Yg\},
$$
being a consequence of the properties of the Schouten bracket: $[[[Y,\eta],g],f]=[Y,[[\eta,g],f]]-[[\eta,g],Yf]-[[\eta,Yg],f]$.

Now choose any  system of linear coordinates $\{x_i\}$ on $\g$ and consider a map $\phi:\g\to\g$, $x_i\mapsto x_i+t X_i(x)$, where $X=X_i(x) \frac{\d }{\d x_i}$ (cf. \cite[Prop. 4.1]{dufourZung}). Let $\psi:\g\to\g$, $x_i\mapsto x_i+t \psi^{(1)}_i(x)+t^2 \psi^{(2)}_i(x)+\cdots$ be the formal (in $t$) inverse map. Note that by the Taylor formula we have $\psi^{(1)}_i=-X_i$ and $\deg\psi^{(j)}=j+1$.

Calculate the first terms of expansion in $t$ of the Poisson bracket $\{x_i,x_j\}_{\phi_*(\pi_1+t\pi_2)}=$ \linebreak $(\phi^{-1})^*\{\phi^*x_i,\phi^*x_j\}_{\pi_1+t\pi_2}$ corresponding to the Poisson bivector $\phi_*(\pi_1+t\pi_2)$:
\begin{align*}
\psi^*\{\phi^*x_i,\phi^*x_j\}_{\pi_1+t\pi_2}=\psi^*\{x_i+tX_i,x_j+tX_j\}_{\pi_1+t\pi_2}=\\
\psi^*(\{x_i,x_j\}_{\pi_1}+
t(\{x_i,x_j\}_{\pi_2}+\{X_i,x_j\}_{\pi_1}+\{x_i,X_j\}_{\pi_1})
+o(t))=\\
\{x_i,x_j\}_{\pi_1}+t(\{x_i,x_j\}_{\pi_2}+\{X_i,x_j\}_{\pi_1}+\{x_i,X_j\}_{\pi_1}-
\frac{\d \{x_i,x_j\}_{\pi_1}}{\d x_s}X_s)
+o(t)=\\
\{x_i,x_j\}_{\pi_1}+t(\{x_i,x_j\}_{\pi_2}+\{Xx_i,x_j\}_{\pi_1}+\{x_i,Xx_j\}_{\pi_1}-
X\{x_i,x_j\}_{\pi_1})
+o(t)\\=\{x_i,x_j\}_{\pi_1}+o(t);
\end{align*}
here we used the Taylor formula and the fact that $\{x_i,x_j\}_{\pi_1}$ is a linear function.

Now we notice that $\phi_*(\pi_1+t\pi_2)=\pi_1+t^2\pi_3+o(t^2)$, where $\pi_3$ is a cubic bivector and $o(t^2)$ is a formal  bivector of degree $> 3$. The fact that $[\phi_*(\pi_1+t\pi_2),\phi_*(\pi_1+t\pi_2)]=0$ implies that $[\pi_1,\pi_3]=0$ and by the triviality of $H^2(\g,S^3\g)$ we conclude that there exists a cubic vector field $X^{(2)}$ such that $[X^{(2)},\pi_1]=\pi_3$ (the proof of this fact is analogous to that of  Lemma \ref{lem}, cf. also \cite[Lemma 4.5]{dufourZung}).

We further mimic the considerations above showing that $\phi^{(2)}_*(\pi_1+t^2\pi_3+o(t^2))=\pi_1+t^3\pi_4+o(t^3)$, where $\phi^{(2)}:\g\to\g$ is given by $x_i\mapsto x_i+t^2X^{(2)}_i$, $\pi_4$ is a quartic bivector and $o(t^3)$ is a formal  bivector of degree $> 4$.

Obviously, $\phi^{(2)}\circ\phi=x_i+tX_i(x)+t^2X_i^{(2)}(x)+o(t^2)$ and $(\phi^{(2)}\circ\phi)_*(\pi_1+t\pi_2)=\pi_1+t^3\pi_4+o(t^3)$, where $X_i(x)$ and $X_i^{(2)}(x)$ are polynomial of required degree. Continuing this process  we build by induction the formal map $\Phi(t)$ thus proving Items 1, 2.

To prove Item 3 notice that $\Phi(t)^*\{f,g\}_{\pi_1}=\{\Phi(t)^*f,\Phi(t)^*g\}_{\pi_1+t\pi_2}$, where $\{f,g\}_\eta$ stands for the Poisson bracket corresponding to the Poisson bivector $\eta$, hence $\{f,\cdot\}_{\pi_1}=0$ implies $\{\Phi(t)^*f,\cdot\}_{\pi_1+t\pi_2}=0$.

Formula (\ref{ca})  follows from the Taylor formula\footnote{Alternatively, one observes that the operator $[X,\cdot]$ is a derivation with respect to the Schouten bracket. Consequently, $[\pi_1,f]=0$ implies $[[X,\pi_1],f]+[\pi_1,Xf]=0$, i.e. $\pi_2(f)+\pi_1(Xf)=0$.}.
\qed

\begin{rema}\label{r-mat}\rm
We conclude this section by noticing that the $r$-matrix type quadratic bivectors $r=r^{ij}\pi_1(x_i)\wedge
\pi_1(x_j)$ mentioned in the introduction are of coboundary type. Indeed, one can take $X:=r^{ij}x_i
\pi_1(x_j)$.
\end{rema}


\section{Centrally linearizable Poisson pencils: general results}
\label{adm}

Let $\s$ be a semisimple Lie algebra over a field $\K=\R,\C$ and let $\g=\s\oplus\K$ be its one-dimensional trivial central extension.
 Choose a system of linear coordinates $x_1,\ldots,x_n$ on $\s^*$ and an element $a\in\K^* \subset\g^*=\s^*\oplus\K^*$, $a\not=0$.  Then $x_0,x_1,\ldots,x_n$ is a system of linear coordinates on  $\g^* $, where $x_0$ is given by  $x_0(a)=1,x_0|_{\s^*}=0$. Note that the coordinate $x_0$ and the vector field $\frac{\d }{\d x_0}$ do not depend on the choice of the coordinates on $\s^*$.

\begin{theo}\label{th1}
 Any CL Poisson bivector $\pi_2$ on $\g^*$ is the Lie derivative
\begin{equation}\label{form}
\pi_2=[\widetilde{X},\pi]=\L_{\widetilde{X}}\pi_1=\pi-\pi_1(q)\wedge \frac{\d }{\d x_0}+x_0\pi_1
\end{equation}
of \label{th1} the canonical Lie--Poisson bivector $\pi_1$ on $\g^*$ along the quadratic vector field $\widetilde{X}$  of the form
$$
\widetilde{X}=X+q\frac{\d }{\d x_0}-x_0E,
$$
where $X$ and $q$ are respectively  a quadratic vector field and quadratic homogeneous polynomial depending only on the coordinates $x_1,\ldots,x_n$, $E=\sum_{i=1}^nx_i\frac{\d }{\d x_i}$ is the Euler vector field on $\s^*$ (treated as a vector field on  $\g^*$), $\pi=[X,\pi_1]=\L_X\pi_1$, and the following identities are satisfied
\begin{equation}\label{eqbasic}
[\pi,\pi]=2\pi_1(q)\wedge\pi_1
\end{equation}
and
\begin{equation}\label{eqbasic2}
[\pi,\pi_1(q)]=0.
\end{equation}
Here $\pi_1(q)=[\pi_1,q]$ stands for the hamiltonian vector field corresponding to the hamiltonian function $q$ and both the identities can be treated as those on $\s^*$.
\end{theo}

\noindent Let $V$ be a vector space, $a\in V$. We denote by $T_a$ the infinitesimal generator of the one-parametric group of transformations $x\mapsto x+ta$. If $R$ is a quadratic polyvector on $V$, then $[T_a,R]$ is a linear polyvector called the {\em linearization }of $R$ at $a$ and if $R$ is a linear polyvector, then $[T_a,R]=R|_a$ (a constant polyvector). Note that $[T_a,[T_a,R]]=R|_a$ for a quadratic $R$.

Let $\pi_2$ be an CL quadratic Poisson bivector on $\g$. Then $\pi_2|_a=0$, $[T_a,\pi_2]=\pi_1$, $[T_a,\pi_1]=\pi_1|_a=0$, where $T_a=\frac{\d }{\d x_0}$. By Lemma \ref{lem} there exists a quadratic vector field $\widetilde{X}$ sucht that $\pi_2=[\widetilde{X},\pi_1]$. Thus $\pi_1=[T_a,[\widetilde{X},\pi_1]]=
[[T_a,\widetilde{X}],\pi_1]+[\widetilde{X},[T_a,\pi_1]]=[[T_a,\widetilde{X}],\pi_1]$.

In particular, $[T_a,\widetilde{X}]=-E+Y_1$, where $E$ is the Euler vector field on $\s^*$ (treated as a vector field on $\g^*$), $[E,\pi_1]=-\pi_1$, and $Y_1$ is any linear vector field with $[Y_1,\pi_1]=0$. Then $Y_1$ is a sum of arbitrary fundamental vector field $X_\xi,\xi\in\g$, of the coadjoint action and a vector field of the form $fT_a$, where the function $f$ is a linear Casimir function of the bivector $\pi_1$, $f=\al x_0$, $\al\in\K$, cf. Remark \ref{rema}.

Let $\widetilde{X}=X_0+x_0X_1+(x_0^2/2)X_2$, where the coefficients of $X_i$ are independent of $x_0$. Then the linearization $[T_a,\widetilde{X}]=X_1+x_0X_2$ is of the form $-E+X_\xi+\al x_0T_a$ and we conclude that $\widetilde{X}$ itself is of the form $\widetilde{X}=X_0-x_0E+x_0X_\xi+\al (x_0^2/2)T_a$. Note that the term $Z:=x_0X_\xi+\al (x_0^2/2)T_a$ can be omitted since $[Z,\pi_1]=0$ (cf. Lemma \ref{lem}(2), Remark \ref{rema}). Finally, we can write
$$
\widetilde{X}=X+qT_a-x_0E,
$$
where $X$ is a quadratic vector field on $\s^*$ understood as that on $\g^*$ (i.e. $X=X_i(x_1,\ldots,x_n)\frac{\d }{\d x_i}$, $i>0$) and $q=q(x_1,\ldots,x_n)$ is a homogeneous quadratic polynomial independent of $x_0$. Consequently,
\begin{equation*}\label{eee}
\pi_2=[\widetilde{X},\pi_1]=\pi-\pi_1(q)\wedge T_a+x_0\pi_1,
\end{equation*}
where $\pi=[X,\pi_1]$ is  a quadratic bivector on $\s^*$ understood as that on $\g^*$ (i.e. $\pi=\pi_{ij}(x_1,\ldots,x_n)\frac{\d }{\d x_i}\wedge \frac{\d }{\d x_j}$, $i,j>0$) and $\pi_1(q)=[\pi_1,q]$ is the hamiltonian vector field with the hamiltonian $q$. Since the vector fields $\pi_1(q)$ and $T_a$ commute, we have $[\pi_1(q)\wedge T_a,\pi_1(q)\wedge T_a]=0$. Also $[x_0\pi_1,x_0\pi_1]=0$, $[\pi,T_a]=0$, $[\pi_1,\pi_1(q)]=0$,  and $[\pi,\pi_1]=0$, whence
\begin{align*}
[\pi_2,\pi_2]=[\pi,\pi]-2[\pi,\pi_1(q)\wedge T_a]+2x_0[\pi,\pi_1]-2[\pi_1(q)\wedge T_a,x_0\pi_1]=\\
[\pi,\pi]-2[\pi,\pi_1(q)]\wedge T_a-2\pi_1(q)\wedge\pi_1=0.
\end{align*}
The second term is the only  containing $T_a$, hence $[\pi_2,\pi_2]=0$ if and only if
\begin{equation*}
[\pi,\pi]=2\pi_1(q)\wedge\pi_1,[\pi,\pi_1(q)]=0.
\end{equation*}
\qed

\begin{theo}\label{th2} If dimension of the generic coadjoint orbit  of the Lie algebra $\s$ is greater or equal 6,  identity (\ref{eqbasic2}) follows from (\ref{eqbasic}).
\end{theo}

\noindent The graded Jacobi identity for the Schouten bracket applied to bivectors $a,b,c$  looks standardly:
$$
[a,[b,c]]+[c,[a,b]]+[b,[c,a]]=0.
$$
Putting $a=b=c=\pi$ we get $[\pi,[\pi,\pi]]=0$. Applying the operator $[\pi,\cdot]$ to identity (\ref{eqbasic}) and using the Leibnitz identity for the Schouten bracket and the fact that $[\pi,\pi_1]=0$ we obtain
$$
[\pi,\pi_1(q)]\wedge\pi_1=0.
$$
Now fix a point $x$ in a generic symplectic leaf of $\pi_1$ and use the Darboux basis for  $\pi_1|_x$:
$$
\pi_1|_x=e_1\wedge e_2+e_3\wedge e_4 +\cdots+e_{2n-1}\wedge e_{2n}.
$$
It is easy to see that under the restriction  $2n\ge 6$ the equation $a\wedge\pi_1(x)=0$, where $a$ is an unknown bivector, has only the trivial solution. Indeed, let $e_1,\ldots,e_{2n},e_{2n+1},\ldots,e_m$ be a basis of $T_x\s^*$ and let $a=\sum_{i< j}a_{ij}e_i\wedge e_j$. Then the equality $a\wedge\pi_1(x)=0$ is equivalent to the system of linear equations
\begin{align*}
a_{12}+a_{34}=0, & a_{12}+a_{56}=0, \ldots, & a_{12}+a_{2n-1,2n}=0,\\
 & a_{34}+a_{56}=0, \ldots, &  a_{34}+a_{2n-1,2n}=0,\\
 &    & \vdots\\
 &     & a_{2n-3,2n-2}+a_{2n-1,2n}=0,\\
 & & a_{ij}=0,
\end{align*}
where $(i,j)$ runs through all the pairs of indices not appearing in the equations above. This system has only the trivial solution.
\qed

\begin{theo}\label{th3} Let $q$ be a quadratic polynomial on $\s^*$ corresponding to a CL quadratic bivector $\pi_2$ by Theorem \ref{th1}. Then
there exists a cubic homogeneous polynomial   $m$ on $\s^*$  such that
\begin{equation}\label{casi2}
[X,\pi_1(q)]=\pi_1(m);
\end{equation}
here we treat $\pi_1$ as a bivector on $\s^*$.
\end{theo}

\noindent  Identity (\ref{eqbasic2}) can be rewritten as
\begin{equation}\label{cubpol}
[\pi_1,[X,\pi_1(q)]]=0
\end{equation}
and can be treated as an identity on $\s^*$. The vector field $[X,\pi_1(q)]$ is cubic and, since $H^1(\s,S^3(\s))=0$, (\ref{cubpol}) holds if and only if there exists a cubic polynomial $m$ on $\s^*$ such that
$$
[X,\pi_1(q)]=\pi_1(m)
$$
(we used the fact that the Poisson cohomology with respect to  $\pi_1$ restricted to the cubic vector fields is naturally identified with the Chevalley--Eilenberg cohomology $H^1(\s,S^3\s)=0$,  cf. \cite[Lemma 4.5]{dufourZung}, Remark \ref{rema} and the proof of Lemma \ref{lem}).
\qed

\begin{theo}\label{th4a}
Let $\pi_2$ be a CL  quadratic Poisson bivector of the form (\ref{form}). Then for any Casimir function $f$ of the bivector $\pi_1$, $\pi_1(f)=0$, there exists a Magri--Lenard chain
$$
f_0:=f,f_1,\ldots
$$
(see Definition \ref{defML}). Taking $f=x_0$, we get a Magri--Lenard chain
$$
x_0,Xx_0=q,f_2,\ldots,
 $$
 which in particular means that the quadratic polynomial $q$ is contained in the canonical involutive family of functions related to the Poisson pencil $\{\la_1\pi_1+\la_2\pi_2\}$ (see Definition \ref{canoo}).
\end{theo}

\noindent The proof follows from Lemma \ref{lem22}(3). \qed

\begin{theo}\label{th4}
Let $\pi_2$ be a CL  quadratic Poisson bivector of the form (\ref{form}) and
let $r=r_0+x_0r_1+\cdots+x_0^kr_k$ be a polynomial of degree $k$ in $x_0$, where $r_0,\ldots,r_k$ are functions of $x_1,\ldots,x_n$. Then
\begin{enumerate}\item
$r$ is a Casimir function of $\pi_2$ if and only if
\begin{align}
\pi_1(q)r_0=\pi_1(q)r_1=\cdots=\pi_1(q)r_k=0,\label{alge2}\\
\pi(r_0)-r_1\pi_1(q)=0,\label{alge}\\
\pi(r_1)-2r_2\pi_1(q)+\pi_1(r_0)=0,\label{alge1}\\
\vdots \nonumber\\
\pi(r_{k-1})-kr_k\pi_1(q)+\pi_1(r_{k-2})=0,\nonumber\\
\pi(r_k)+\pi_1(r_{k-1})=0,\label{a3a}\\
\pi_1(r_k)=0\label{a3},
\end{align}
where the above identities can be interpreted as those on $\s^*$;
\item if equality (\ref{a3}) is satisfied for $r_k$,  then equality (\ref{a3a}) is satisfied for $r_{k-1}=Xr_k+ r'_{k-1}$, where $r'_{k-1}$ is a Casimir function for $\pi_1$; moreover, $\pi_1(q)r_{k-1}=0$.
    \end{enumerate}
\end{theo}

\noindent  Indeed, by (\ref{form}) we have
\begin{align*}
0=\pi_2(r)=\pi(r_0)+x_0\pi(r_1)+\cdots+x_0^k\pi(r_k)\\
+(\pi_1(q)r_0+x_0\pi_1(q)r_1+\cdots+x_0^k\pi_1(q)r_k)\frac{\d }{\d x_0}\\
-(r_1+2x_0r_2+\cdots+kx_0^{k-1}r_k)\pi_1(q)\\
+x_0\pi_1(r_0)+x_0^2\pi_1(r_1)+\cdots+x_0^{k+1}\pi_1(r_k),
\end{align*}
which implies Item 1.

Applying the operator $[X,\cdot]$ to the equality $\pi_1(r_k)=[\pi_1,r_k]=0$, we get $\pi(r_k)+\pi_1(Xr_k)=0$ and we can put $r_{k-1}:=Xr_k$.  From the obvious equalities  $[\pi_1(m),r_k]=0$, $[\pi_1(q),r_k]=0$ and equality (\ref{casi2}) we get
\begin{equation}
\label{A0}
0=[\pi_1(m),r_k]=[[X,\pi_1(q)],r_k]=\pi_1(q)(Xr_k)=\pi_1(q)r_{k-1},
\end{equation}
which finishes the proof of Item 2.
\qed

\begin{coro}\label{cor} Assume that equalities of Theorem \ref{th4}(1) hold. Then
$\pi_1(r_i)r_j=0=\pi(r_i)r_j$ for any $i,j$.
\end{coro}

\noindent  Indeed, $\pi_1(r_i)r_j=-\pi(r_{i+1})r_j=\pi(r_{j})r_{i+1}=-\pi_1(r_{j-1})r_{i+1}$. Continuing this process we will finally come to possibility to use  equality (\ref{a3}). \qed

\begin{theo}\label{th6}
Let $\pi_2$ be a CL  quadratic Poisson bivector of the form (\ref{form}) and
let $r(x_0):=r_0+x_0r_1+\cdots+x_0^kr_k$ be a  Casimir function of the bivector $\pi_2$ being a polynomial of  degree $k$ in $x_0$, where $r_0,\ldots,r_k$ are functions of $x_1,\ldots,x_n$. Then
\begin{enumerate}\item
the function
$$
C(\la)=r(x_0+\la)=r(x_0)+r'(x_0)\la+\cdots+\frac{r^{(k)}(x_0)}{k!}\la^k
$$
is a Casimir function of the Poisson bivector $\pi_2+\la\pi_1$;
\item the functions $f_0:=r(x_0),f_1:=r'(x_0),\ldots,f_k:=\frac{r^{(k)}(x_0)}{k!}$
satisfy the relations
$$
\pi_2(f_0)=0,\pi_1(f_0)+\pi_2(f_1)=0,\ldots,\pi_1(f_{k-1})+\pi_2(f_k)=0,\pi_1(f_k)=0
$$
in particular $f_k,f_{k-1},\ldots,f_0,0,0,\ldots$ is a Magri--Lenard chain in the sense of Definition  \ref{defML}.
\end{enumerate}
\end{theo}

\noindent It is well-known that, given a quadratic Poisson bivector $\pi$ on a vector space $V$ vanishing at some point $a\in V$, $a\not=0$, the linearization $\pi_a$ of $\pi$ at $a$ is a linear Poisson bivector compatible with $\pi$ and, moreover, the space $C(\pi+\la\pi_a)$ of Casimir functions of the Poisson bivector $(\pi+\la\pi_a)(x)=\pi(x+\la a)$ is equal to $\{f(x+\la a)\mid f\in C(\pi)\}$. In particular, a CL Poisson bivector $\pi_2$ is compatible with the Lie--Poisson bivector $\pi_1$ and the Casimir functions of the Poisson bivector $\pi_2+\la\pi_1$ are the shifts of the Casimir functions of $\pi_2$ in direction of the element $a$. In this context the proof Item (1) is a  consequence of the fact that $x_0(a)=1$ and $(\pi_2+\la\pi_1)(x_0,\ldots,x_n)=\pi_2(x_0+\la,x_1,\ldots,x_n)$.

Item (2) follows from Item (1).
\qed

\begin{rema}\label{rre}\rm
By Theorem \ref{th4a}, given any Casimir function $f$ of the bivector $\pi_1$, there exists a Magri--Lenard chain $f_0=f,f_1,\ldots$, i.e. a formal Casimir function $f_0+\la f_1+\cdots$ of the bivector $\pi_1+\la\pi_2$. If $f$ is a homogeneous polynomial of degree $k$, then Lemma \ref{lem22} shows that the functions $f_i$ could be chosen homogeneous polynomials of degree $k+i$. We say that such a Magri--Lenard chain \emph{is of degree} $k$. If $f_0,f_1,\ldots$ and $g_0,g_1,\ldots$ are two Magri--Lenard chains of degrees $k$ and $m$ respectively, $m\ge k$, then obviously, $f_0,\ldots,f_{m-1},f_m+g_0,f_{m+1}+g_1,\ldots$ is a Magri--Lenard chain of degree $k$, hence the initial data $f_0$ defines the chain up to the addition of other chains of corresponding degrees. Using this observation one can search for the homogeneous polynomial Casimir functions of the bivector $\pi_2$ starting from a fixed Magri--Lenard chain and adding linear combinations of other Magri--Lenard chains, see Remark \ref{casi22}.
\end{rema}

\begin{rema}\label{kro}\rm
This remark concerns the Jordan--Kronecker decomposition of the CL Poisson pencil $\{\la_1\pi_1+\la_2\pi_2\}$ \cite{gz4}, \cite{pSymmetries}, \cite{bolsSing}. Assume for simplicity the complex case. It follows from Theorem \ref{th4a} that this decomposition does not contain Jordan blocks with the constant eigenvalue equal to zero (otherwise   the corank of the bivector $\pi_1$ at a generic point would be greater than that of $\pi_1+\la\pi_2$ for generic $\la\in\C$ and for some Casimir functions $f$ of the bivector $\pi_1$  it would be impossible to find a Magri--Lenard chain starting from $f$, cf. Remark \ref{remmm}). On the other hand, the fact that all the bivectors $\pi_1+\la\pi_2$, $\la\in\C\cup\{\infty\}$, $\la\not=0$, are diffeomorphic to each other (by means of the  translation along $a$) implies that there are no Jordan blocks with constant eigenvalues at all (otherwise for some $\la$ corank of $\pi_1+\la\pi_2$ would jump up). These two observations allow to come to the following conclusion: in order to prove the kroneckerity of the pencil  (recall that this means the absence of the Jordan blocks at a generic point and guarantees the completeness of the canonical involutive family of functions) it is enough to show that the corank of the singular set $\Sing\pi_2$ of the Poisson bivector $\pi_2$ (or any other quadratic bivector of the pencil, as they all are diffeomorphic), i.e. the union of all its symplectic leaves of nonmaximal dimension, is of codimension at leat two.

Consider the CL Poisson bivector (\ref{tauteq}) from Theorem \ref{tauto} below and assume that the dimension of generic coadjoint orbits of the Lie algebra $\s$ is greater than 2. Then $\Sing\pi_2$ contains the hyperplane $x_0-x=0$ and the condition $\codim\Sing\pi_2\ge 2$ is not satisfied. This in particular means the presence of Jordan blocks with eigenvalues depending on the point and noncompleteness of the canonical involutive family of functions.
\end{rema}


\section{First examples}
\label{taut}

We have two limiting possibilities for solving the basic identiies \ref{eqbasic}, \ref{eqbasic2}.

\begin{exa}\label{sklya}
\rm {\em (Sklyanin bracket)}
First is when  $X$ and $\pi$ are zero. Then (\ref{eqbasic}) is equivalent to $\pi_1(q)\wedge\pi_1=0$. The case when $\pi_1(q)=0$ is rather noninteresting since it produces a Poisson pencil of bivectors proportional to $\pi_1$. An example of such a situation with  $\pi_1(q)\not=0$ is as follows. Let $\s=\so(3)$, $\pi_1=x_1 \frac{\d }{\d x_2}\wedge \frac{\d }{\d x_3}+x_2 \frac{\d }{\d x_3}\wedge \frac{\d }{\d x_1}+x_3 \frac{\d }{\d x_1}\wedge \frac{\d }{\d x_2}$, $q=\frac{1}{2}(J_1x_1^2+J_2x_2^2+J_3x_3^2)$, where $J_i$ are arbitrary numbers. Then $\pi_1(q)\wedge\pi_1=0$ and $\pi_2=x_0\pi_1+x_2x_3(J_2-J_3)\frac{\d }{\d x_0}\wedge \frac{\d }{\d x_1}+x_1x_2(J_1-J_2)\frac{\d }{\d x_0}\wedge \frac{\d }{\d x_3}+x_3x_1(J_3-J_1)\frac{\d }{\d x_0}\wedge \frac{\d }{\d x_2}$ is the Sklyanin quadratic Poisson bivector on $\g$ \cite{sklyanin}, \cite[Ex. 7]{GMP}.

Note however that such a situation is rather exceptional  due to the low dimensionality  of $\s$. In general quadratic functions $q$ with $\pi_1(q)\not=0$ and $\pi_1(q)\wedge\pi_1=0$ do not exist.
\end{exa}

\begin{exa}\label{sklya}
\rm
Second possibility  is when $q$ is a Casimir function of $\pi_1$. Then $\widetilde{X}=X-x_0E,\pi_2=\pi+x_0\pi_1$, where $\pi=[X,\pi_1]$, and (\ref{eqbasic}) is equivalent to $[\pi,\pi]=0$. In particular, one can take $\pi$ to be one of the $r$-matrix type quadratic bivectors $r^{ij}\pi_1(x_i)\wedge
\pi_1(x_j)$, $i,j>0$, mentioned in Introduction (see also Remark \ref{r-mat}).
Note that here any  Casimir function of $\pi_1$ is a Casimir function of $\pi_2$ and consequently of any bivector $\pi_1+\la\pi_2$.
\end{exa}

Next theorem gives a series of examples which we call \emph{the tautological solutions of the  basic identities.}

\begin{theo}\label{tauto} Let $\s$ be a Lie algebra, $E$ the Euler vector field on $\s^*$, $x\in\s$, a linear function on $\s^*$. Then the quadratic vector field $X:=xE$ along with the quadratic polynomial $q:=x^2/2$ satisfy identities (\ref{eqbasic}), (\ref{eqbasic2}). Consequently the  bivector
\begin{equation}\label{tauteq}
\pi_2=[X,\pi_1]-\pi_1(q)\wedge \frac{\d }{\d x_0}+x_0\pi_1=(x_0-x)\pi_1-\pi_1(x)\wedge(E+x \frac{\d }{\d x_0})
\end{equation}
is a CL  Poisson bivector on $\g^*=(\s\oplus\K)^*$ when $s$ is semisimple.
\end{theo}

\noindent Using the facts that  $[E,\pi_1]=-\pi_1$, $[E,x\pi_1]=0$, $[E,x]=x$ and $[E,Y]=0$ for any linear vector field $Y$, we get $[X,\pi_1]=[xE,\pi_1]=-[\pi_1,xE]=-\pi_1(x)\wedge E-x[\pi_1,E]=-\pi_1(x)\wedge E-x\pi_1$,
$[\pi_1(x)\wedge E,\pi_1(x)\wedge E]=0$, $[\pi_1(x)\wedge E,x\pi_1]=0$, $[x\pi_1,x\pi_1]=2x\pi_1(x)\wedge\pi_1$.
Thus
$$
[[X,\pi_1],[X,\pi_1]]=[-\pi_1(x)\wedge E-x\pi_1,-\pi_1(x)\wedge E-x\pi_1]=[x\pi_1,x\pi_1]=2x\pi_1(x)\wedge\pi_1=2\pi_1(x^2/2)\wedge \pi_1.
$$
Moreover,
$$
[[X,\pi_1],\pi_1(x^2/2)]=[-\pi_1(x)\wedge E-x\pi_1,x\pi_1(x)]=-[x\pi_1,x\pi_1(x)]=0.
$$
\qed

\section{Particular irreducible modules in the space of quadratic vector fields and quadratic polynomials on $\sl(3,\C)^*$}
\label{s10}

Let $\s=\sl(3,\C)$ and let $\om_1,\om_2$ be the system of fundamental weights of $\s$. We denote by $R(\al)$ the representation of $\s$ with the highest weight $\al$. In particular,  $V=R(\om_1)$ will stand for the natural $3$-dimensional linear representation of $\s$, $V^*=R(\om_2)$ for its dual, $\s=\s^*=R(\om_1+\om_2)$ for the adjoint representation. The space $\Vect^2(\s^*)$ of quadratic vector fields on $\s^*$ is naturally identified with $S^2(\s)\otimes\s^*$.

According to \cite[Table 5]{vinbergOniLieGrAlgGr} the space $S^2(\s)$ of quadratic polynomials on $\s^*$ is decomposed to the following irreducible components: $S^2(\s)= R(2\om_1+2\om_2)+\s+1$, where $1=R(0)$ stands for the trivial 1-dimensional representation. It can be shown (LiE program) that  the component $R(2\om_1+2\om_2)\otimes\s$  of the module $\Vect^2(\s^*)$ contains the submodules $W:=R(3\om_1)=S^3(V)$ and $W^*=R(3\om_2)=S^3(V^*)$. An explicit realization of the module $W^*$ as a submodule of $R(2\om_1+2\om_2)\otimes \s$ is indicated in Appendix and from now on we mean by $W^*$ this realization\footnote{\label{secc} There are also submodules in $\s\otimes\s$ isomorphic to $R(3\om_1)$ and $R(3\om_2)$, which however do not lead to the solution of basic identities and we do not consider them.}.

Moreover,
\begin{equation}\label{dec}
S^2(W^*)=R(2\om_1+2\om_2)+R(6\om_2)
\end{equation}
 and $S^2(W)=R(2\om_1+2\om_2)+R(6\om_1)$. Dimensions of the spaces $W$ and $R(2\om_1+2\om_2)$ are 10 and 27 respectively  \cite[Table 5]{vinbergOniLieGrAlgGr}. All the weight spaces of $W$ and $W^*$ are of multiplicity one and the following figure shows the weight spaces and multiplicities of the module $R(2\om_1+2\om_2)$ in the positive Weyl chamber.
$$
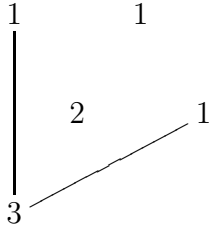
\xymatrixcolsep{1pc}\xymatrix{
& & & & 1   &   & 1 & &\\
 & & &  & & 2 & & 1&\\
& &  &   &  3 \ar@{-}[uu] \ar@{-}[rrru] & & & & \\
}
$$
The explicit realization of the module $R(2\om_1+2\om_2)$, which will be denoted by $U$, in the space of quadratic polynomials on $\sl^*(3)$ can be described by the following generators:
$$
\xymatrixcolsep{-1pc}\xymatrix{
& & x_{23}^2& & x_{23}x_{13}   &   & x_{13}^2 & &\\
 &x_{21}x_{23} & & P_{23}q_{21},P_{23}q_{13} & & P_{13}q_{23},P_{13}q_{12} & & x_{12}x_{13}&\\
x_{21}^2& & P_{21}q_{23},P_{21}q_{31} &   &  q_{12},q_{13},q_{23}& & P_{12}q_{13},P_{12}q_{32}& & x_{12}^2\\
& x_{21}x_{31}&  &P_{31}q_{12},P_{31}q_{32}& & P_{32}q_{12},P_{32}q_{31}& & x_{12}x_{32} & \\
& &x_{31}^2 & &  x_{31}x_{32}& &x_{32}^2, & &}
$$
where
$$
q_{ij}=x_{ij}x_{ji}-\frac12(x_{ii}-x_{jj})^2,
P_{ij}q_{jk}=x_{ik}x_{kj}-x_{ij}(x_{jj}-x_{kk}),
P_{ij}q_{ik}=-x_{ik}x_{kj}+x_{ij}(x_{ii}-x_{kk}).
$$

\section{10-parametric family of CL quadratic bivectors on $\gl(3)^*$}
\label{ser}

Let $P$ be the canonical Lie--Poisson bivector on $\sl(3,\K)^*$, $\K=\R,\C$.

\begin{theo}\label{mainId}
Let $X_0,\ldots,X_{9}$ be the generators of the weight spaces of 10-dimensional irreducible representation $W^*$ described in Section \ref{s10} and Appendix. Then for any $b=(b_0,\ldots,b_{9})\in\K^{10}$ there exists a quadratic polynomial $Q_b\in U$ such that
\begin{equation}\label{id1}
[[X_b,P],[X_b,P]]=2P(Q_b)\wedge P,
\end{equation}
where $X_b:=b_0X_0+\cdots+b_{9}X_{9}$. The explicit formula for $Q_b$ is given in Appendix\footnote{Although we defined the modules $W^*$ and $U$ over $\C$, all the coefficients in the formulas in Appendix are real, hence the theorem is valid over $\R$ as well.}.
\end{theo}

\noindent Notice that the correspondence $(X,Y)\mapsto[[X,P],[Y,P]]$ defines a symmetric bilinear $SL(3)$-equivariant map from $W^*\times W^*$ to the space $S^3(\s)\otimes\bigwedge^3(\s^*)$ of cubic trivectors on $\s^*$, hence can be factorized to a linear $SL(3)$-equivariant map $\phi:S^2(W^*)\to S^3(\s)\otimes\bigwedge^3(\s^*)$. By (\ref{dec}) the image of $\phi$ can be isomorphic to one of the modules ${0},U,R(6\om_2),U+R(6\om_2)$ and we will check that in fact  $\phi(S^2(W^*))\cong U$.

First check that (\ref{id1}) holds at least for one $b$ with $Q_b\in U$, $Q_b\not=0$. And indeed direct calculations show that
\begin{equation}\label{che}
[[X_0,P],[X_0,P]]=2P(24q_0)\wedge P,
\end{equation}
where $X_0$ is a generator of the zero weight space in $W^*$ given explicitly in Appendix
and $q_0=-\frac1{4}(q_{12}+q_{23}+q_{13})$.

Now  (\ref{che}) implies that $\phi(S^2(W^*))$ is nontrivial and the following arguments can be used to show that it is isomorphic to $U$. In the space $S^2(W^*)=S^2(S^3V^*)$ of quadratic polynomials on $W=S^3V$ there is a subspace consisting of those vanishing on the image of the Veronese embedding $v_3(V^*) \subset S^3V^*$, $v_3(x)=x\odot x\odot x$. This is exactly $R(2\om_1+2\om_2)$ (the complement subspace $R(6\om_2)=S^6V^*$ of 6-th order polynomials on $V^*$ is a linear span of the image of the composition of two Veronese embeddings $v_2\circ v_3(V^*)$, $v_2:W^*\to S^2(W^*)$, $v_2(X)=X\odot X$, cf. (\ref{dec}) and \cite[\S 13]{fultonharris}). In particular, in order to show that $\phi(S^2(W^*))\cong R(2\om_1+2\om_2)$ it is enough to ensure that $\phi(X)=0$ for $X\in v_3(V^*)$. In turn, it is enough to check the last fact for the highest weight vector $X$ of the representation $W^*$ (since the image of the Veronese embedding $v_3$ coincides with the $SL(3)$-orbit of the highest weight vector). Indeed, direct check shows that $[[X,P],[X,P]]=0$ for the  highest weight vector $X=X_8$ (see Appendix) of the representation $W^*$.

Now, since the module $U':=\{\pi_1(q)\wedge\pi_1\mid q\in U\}\subset S^3(\s)\otimes\bigwedge^3(\s^*)$ is isomorphic to $U$ and $U'\cap \phi(S^2(W^*))\not =\{0\}$, we conclude by the equivariance that $Q_b$ is simply $\frac1{2}\mu\circ\phi\circ v_2(X_b)$, where $v_2:W^*\to S^2(W^*)$ is the Veronese embedding and $\mu:\phi(S^2(W^*))\to U$ is the corresponding isomorphism.
\qed

\medskip

 Theorems \ref{mainId}, \ref{th1} and \ref{th2} imply  the following theorem.

\begin{theo}\label{th7} Let $x_0$ stand for any coordinate on $\gl^*(3,\K)=\sl(3,\K)^*\oplus\K^*\cong\sl(3,\K)\oplus\K=\gl(3,\K)$ proportional to the function $\Trace(\cdot)$ and let $X_b$ and $Q_b$ be the quadratic vector field and the quadratic polynomial from Theorem \ref{mainId} treated as objects on $\gl^*(3,\K)$ independent of this coordinate.
Put $\widetilde{X}_b=X_b+Q_b \frac{\d }{\d x_0}-x_0E$. Then $\pi_2(b):=[\widetilde{X}_b,\pi_1]=[X_b,\pi_1]-\pi_1(Q_b)\wedge \frac{\d }{\d x_0}+x_0\pi_1$ is a CL quadratic Poisson bivector on $\gl(3,\K)^*$, where $\pi_1$ is the Lie--Poisson bivector of $\gl(3,\K)^*$\footnote{\label{footnn} The Maple software calculates the Schouten bracket of two bivectors $P=P^{ij} \frac{\d }{\d x^i }\wedge \frac{\d }{\d x^j}$ and $R=R^{ij} \frac{\d }{\d x^i }\wedge \frac{\d }{\d x^j}$ by the formula $[P,R]^{ijk}=2P^{r[i}\frac{\d }{\d x^r}R^{jk]}-2R^{r[k}\frac{\d }{\d x^r}P^{ji]}$ \cite[Form. (2.3)]{Nij1} (consequently $[P,P]^{ijk}=\frac{4}{3}\sum_{c.p. i,j,k}P^{ri}\frac{\d }{\d x^r}P^{jk}$, where the sum is taken over the cyclic permutations of the indices $i,j,k$). In particular, this gives   $[\pi_1(q)\wedge \frac{\d }{\d x_0},x_0\pi_1]=-\frac{2}{3}\pi_1(q)\wedge \pi_1$, cf. the proof of Theorem \ref{th1}. Thus  in Maple, in order that the vector field $\widetilde{X}_b$ generated a CL Poisson bivector it should be of the form $\widetilde{X}_b=X_b-\frac{3}{2}Q_b \frac{\d }{\d x_0}-x_0E$, where $X_b$ and $Q_b$ satisfy (\ref{id1}). Alternatively $\widetilde{X}$ can be of the form from Theorem \ref{th7}, but $X_b$ and $Q_b$ should be related by $[[X_b,P],[X_b,P]]=-\frac{4}{3}P(Q_b)\wedge P$.}.
\end{theo}

\begin{rema}\label{sokk}\rm
The members of the 3-parametric family of quadratic vector fields corresponding to the  Sokolov CL quadratic Poisson bivectors on $\gl(3,\K)^*$ \cite{SokolovEllCM} (see Introduction) modulo inessential quadratic hamiltonian vector fields in our coordinates look as follows: $2X_7-\frac1{2}X_4+\frac1{2}g_2X_5+\frac1{2}g_3X_9$, where $X_i$ are the quadratic vector fields from the module $W^*$ given explicitly in Appendix.
\end{rema}

\begin{rema}\rm
One can easily ``dualize'' the results of this section using the quadratic vector fields from the dual module $W=(W^*)^*$ (on practice this could be done by changing $x_{ij}$ to $x_{ji}$ in all the formulas). We do not consider this case separately, since it gives the same families of functions in involution.
\end{rema}

\section{Casimir functions of the quadratic bivector and the involutive  family of functions of the CL Poisson pencil on $\gl(3)^*$}
\label{cas1}

\begin{theo}\label{cas}
Let $\pi_2=\pi_2(b)$ be the quadratic bivector from Theorem \ref{th7} and let $p:\gl(3,\K)^*=\sl(3,\K)^*\oplus\K^*\to\sl(3)^*$ be the natural projection. Then the following functions are  Casimir functions of $\pi_2$:
\begin{enumerate}\item $K_1=p^*C_3$, where $C_3$ is the cubic Casimir function of the linear Poisson bivector $P$ on $\sl(3)^*$ corresponding to the function $\Trace(x^3)$ under the identification $\sl(3)\cong\sl(3)^*$ by means of the Killing form;
\item $K_2=k+Bx_0$,   where $B:=p^*C_2$, $C_2$ is the quadratic Casimir function of $P$ corresponding to the function $\Trace(x^2)$ and $k=k_b:=p^*X_bC_2$;
\item $K_3:=-3p^*(m-Xq)+3p^*qx_0+x_0^3$, where $m=m_b$ is the cubic homogeneous polynomial corresponding to the vector field $X:=X_b$ and quadratic polynomial $q:=Q_b$ by Theorem \ref{th3}\footnote{In Maple the corresponding function is of the form $K_3:=-\frac{9}{2}p^*(m-Xq)+\frac{9}{2}p^*qx_0+x_0^3$, cf. footnote \ref{footnn}.}. \end{enumerate}
    The explicit formulas for the polynomials $k_b$ and $M_b:=m_b-X_bQ_b$ are given in Appendix.
\end{theo}

\noindent To prove Item 1 observe that $K_1$ does not depend on $x_0$, the coordinate corresponding to the centre, and by Theorem \ref{th4} is a Casimir function of $\pi_2$ if and only if $\pi(C_3)=0$. In order to prove the last equality note that the function $C_3$ is the generator of the trivial 1-dimensional representation $1=R(0)$ in the decomposition $S^3(\s)= R(3\om_1+3\om_2)+R(2\om_1+2\om_2)+\s+R(3\om_1)+R(3\om_2)+1$ (LiE program). Due to the fact that $[P,C_3]=0$ and to the Jacobi identity the equality
\begin{equation}\label{o}
[X_b,C_3]=X_b(C_3)=0
\end{equation}
 implies  $\pi(C_3)=[[X_b,P],C_3]=0$. Thus it is enough to prove (\ref{o}):
on one hand, the function $[X_b,C_3]$ lies in $W^*\otimes 1\cong W^*=R(3\om_2)$, on the other hand there is no such a component in the  space $S^4(\s)$ of quartic polynomials on $\s^*$ (LiE Program).

To prove Item 2 use Theorem \ref{th4} for the cubic polynomial of the form $r=r_0+x_0r_1$. In order that $r$ is a Casimir function the following equalities should be satisfied:
\begin{align}
\pi_1(q)r_0=\pi_1(q)r_1=0,\nonumber\\
\pi(r_0)-r_1\pi_1(q)=0,\label{2q}\\
\pi(r_1)+\pi_1(r_0)=0,\label{2qi}\\
\pi_1(r_1)=0,\nonumber
\end{align}
which we treat as equalities on  $\s^*$.

Put  $r_1=C_2$, $r_0=XC_2$. In view of Theorem \ref{th4}(2) it is enough to show (\ref{2q}).

First check that (\ref{2q}) holds at least for one $b$ with $Q_b\in U$, $Q_b\not=0$. And indeed direct calculations show that (cf. formula (\ref{che}))
\begin{equation}\label{che1}
[[X_0,P],X_0C_2]=24P(C_2q_0),
\end{equation}
where $X_0$ is a generator of the zero weight space in $W^*$ given explicitly in Appendix
and $q_0=-\frac1{4}(q_{12}+q_{23}+q_{13})$\footnote{In Maple this formula looks as $[[X_0,P],X_0C_2]=36P(C_2q_0)$, cf. footnote \ref{footnn}.}.

Now use the argument similar to that from the proof of Theorem \ref{mainId}. Observe that the correspondence $X\mapsto[[X,P],XC_2]=\pi(r_0)$ is the quadratic form of a  bilinear $SL(3)$-equivariant map $(X,Y)\mapsto[[X,P],YC_2]$ from $W^*\times W^*$ to the space $S^4(\s)\otimes\s^*$ of quartic vector fields on $\s^*$, which can be factorized to a linear $SL(3)$-equivariant map $\phi:W^*\otimes W^*=S^2(W^*)\oplus\bigwedge^2(W^*)\to S^4(\s)\otimes \s^*$. We are interested in the image of the restriction of $\phi$ to the diagonal $D:=\{X\otimes X)\mid X\in W^*\}\subset S^2(W^*)$.

By (\ref{dec}) the linear span $\Span(\phi(D))$ of $\phi(D)$ can be isomorphic to one of the modules ${0},U,R(6\om_2),U+R(6\om_2)$ and we have to check that in fact it is isomorphic to $U$.  By (\ref{che1}) $\Span(\phi(D))$ is nontrivial. To show that it is isomorphic to $U$ it is enough to check the equality $[[X,P],XC_2]=0$ on the highest weight vector $X=X_8$ (see Appendix), which can be done directly.

Finally, by the equivariance $[[X_b,P],X_bC_2]=P(C_2Q_b)$ for any $b$.

To prove Item 3 consider a cubic homogeneous polynomial of the form $r=r_0+x_0r_1+x_0^3$ (i.e. $r_2=0,r_3=1$).  By Theorem \ref{th4} it is the  Casimir function if and only if
\begin{align}
\pi_1(q)r_0=\pi_1(q)r_1=0,\nonumber\\
\pi(r_0)-r_1\pi_1(q)=0,\label{2we}\\
\pi(r_1)+\pi_1(r_0)=0,\label{2wer}\\
-3\pi_1(q)+\pi_1(r_1)=0,\nonumber
\end{align}
(identities on $\s^*$).
Let $r_1=3q$. By (\ref{casi2}) we have
$$
\pi(q)=\pi_1(m-Xq),
$$
so one can put $r_0:=-3(m-Xq)$ in order to satisfy (\ref{2wer}). Note that $\pi_1(q)r_0=-\pi_1(r_0)q=\pi(r_1)q=3\pi(q)q=0$ by (\ref{2wer}). It remains to check (\ref{2we}).

 This can be done by direct computer check (as well,  the fact that the function $K_3$ is the Casimir function of $\pi_2$ also can be checked directly).
\qed

\begin{rema}\label{casi22}\rm
Let us use Remark \ref{rre} to explain why the Casimir functions $K_2,K_3$ are of the form given above. Theorem \ref{th4a} shows that, given any Casimir function $f$ of the bivector $\pi_1$, there exists a Margi--Lenard chain $f_0=f,f_1,\ldots$.

Consider $f=B$ and $f=x_0B$. Then  there exist  Magri--Lenard chains   of the form\footnote{Below we abuse notations and treat the vector field $X$ and the functions $m$ and $q$ as objects on  $\gl(3)\cong\gl(3)^*$ independent of $x_0$.} $f_0=B,f_1=XB,f_3,\ldots$ and  $g_0=x_0B,g_1,\ldots$. Try to find $\alpha\in \K$ such that for the Magri--Lenard chain $f_0,f_1+\alpha g_0,f_2+\alpha g_1,\ldots$ the following equality holds: $\pi_2(f_1+\alpha g_0)=0$. We have $\pi_2(f_1+\alpha g_0)=\pi(XB)-\alpha B\pi_1(q)-x_0(\alpha\pi(B)+\pi_1(XB))$, which is zero for $\alpha=1$ by (\ref{2q}) and (\ref{2qi}).

Analogously, take $f=x_0$, $f=x_0^2$ and $f=x_0^3$. Then there exist  Magri--Lenard chains of the form $f_0=x_0,f_1=q,f_2=-(m-Xq),f_3,\ldots$, $g_0=x_0^2, g_1=2x_0q, g_2,\ldots$ and $h_0=x_0^3,h_1,\ldots$ respectively. Try to find $\alpha, \beta \in \K$ such that for the Magri--Lenard chain $f_0,f_1+\alpha g_0,f_2+\alpha g_1+\beta h_0,f_3+\alpha g_2+\beta h_1\ldots$ the following equality holds: $\pi_2(f_2+\alpha g_1+\beta h_0)=0$. We have $\pi_2(f_2+\alpha g_1+\beta h_0)=\pi(-(m-Xq))-2\alpha q\pi_1(q)+x_0(2\alpha\pi(q)+\pi_1(-(m-Xq))+x_0^2(2\alpha\pi_1(q)-3\beta\pi_1(q))$, which is zero for $\alpha=(1/2)$ and $\beta=(1/3)$.

Summarizing this remark we conjecture that in general, given any homogeneous Casimir function $r_k$ of the bivector $\pi_1$ one should be able to construct the functions $r_{k-1},r_{k-2},\ldots,r_0$ satisfying equalities of Theorem \ref{th4}(1) in order to construct a homogeneous Casimir function of the bivector $\pi_2$. To this end one should first construct the Magri--Lenard chains starting from the functions $r_k,x_0r_k,x_0^2r_k,\ldots$ and then look for  a linear combination of these chains with a Casimir funtion of $\pi_2$ at some level. The case of the Casimir function $x_0$, which is important since it incorporates the hamiltonian $q$, is analogous.
\end{rema}

\begin{theo}\label{finv} The functions $x_0,B,K_1,K_2,K_3,p^*Q_b$ Poisson commute with respect to any bivector $\la_1\pi_1+\la_2\pi_2$. In particular,
the functions $Q_b,X_bC_2,m_b-X_bQ_b,C_2,C_3$ form an involutive family of functions with respect to the Lie--Poisson  bivector $P$ on $\sl(3)^*$.
\end{theo}
The proof follows from Theorem \ref{th6} and Remark \ref{invv}. \qed

\begin{rema}\label{comple}
\rm Although we made some preparations for study of completeness of the family of functions from Theorem \ref{finv} (see Remark \ref{kro}), we leave this matter aside. Indeed, it is not easy to prove the condition $\codim \Sing\pi_2\ge 2$: the set $\Sing\pi_2$ consists of the points in $\g^*$, where the differentials of the functions $K_1,K_2,K_3$, depending on 10 parameters, become linearly dependent. However, we conjecture that this family coincides with the canonical involutive family of functions and is complete for any nontrivial value of parameter $b$ and, moreover, the Jordan--Kronecker decomposition of the pencil at a generic point contain three Kronecker blocks of dimensions 1,2 and 3.

Another important question is study of independence/dependence of the canonical commuting families of functions for different values of the parameter $b$.
\end{rema}

\section{Appendix: explicit formulas}
The vector filed $X=b_0X_0+b_1X_1+\cdots+b_9X_9$, where the vector fields $X_i$ are the generators of the irreducible module $W^*$, see Section \ref{s10} (here $y_{13}=x_{11}-x_{33},y_{23}=x_{22}-x_{33}$):
\begin{dgroup*}
\begin{dmath*}
X_0 =((4y_{13}-2y_{23})x_{12}-6x_{13}x_{32})\frac{\d }{\d x_{12}}
+ ((-4y_{13}+2  y_{23})  x_{13}+6  x_{12}  x_{23})\frac{\d }{\d x_{13}}
 + ((2  y_{13}-4  y_{23})  x_{21}+6  x_{23}  x_{31})\frac{\d }{\d x_{21}}
+ ((-2  y_{13}+4  y_{23})  x_{23}-6  x_{13}  x_{21})\frac{\d }{\d x_{23}}
 + ((-2  y_{13}-2  y_{23})  x_{31}-6  x_{32}  x_{21})\frac{\d }{\d x_{31}}
+ ((2  y_{13}+2  y_{23})  x_{32}+6  x_{12}  x_{31})\frac{\d }{\d x_{32}}
 +(-2  y_{13}^2-4  y_{13}  y_{23}+4  y_{23}^2)\frac{\d }{\d y_{13}}
+(-4  y_{13}^2+4  y_{13}  y_{23}+2  y_{23}^2)\frac{\d }{\d y_{23}};
\end{dmath*}
\begin{dmath*}
X_1 =-[P(x_{12}),X_0]=(-3 x_{12}^2)\frac{\d }{\d x_{12}}+(6 x_{13} x_{12})\frac{\d }{\d x_{13}}
 +(-3 x_{12} x_{21}+3 x_{13} x_{31}-3 x_{23} x_{32}+y_{13} y_{23}+y_{23}^2)\frac{\d }{\d x_{21}}+(-2 x_{13} y_{13}+4 x_{13} y_{23})\frac{\d }{\d x_{23}}
 +(3 x_{12} x_{31}-x_{32} y_{13}+5 x_{32} y_{23})\frac{\d }{\d x_{31}}+(-3 x_{32} x_{12})\frac{\d }{\d x_{32}}
 +(2 x_{12} y_{13}+5 x_{12} y_{23}-3 x_{13} x_{32})\frac{\d }{\d y_{13}} + (4 x_{12} y_{13}+x_{12} y_{23}+3 x_{13} x_{32})\frac{\d }{\d y_{23}};
\end{dmath*}
\begin{dmath*}
X_2 =-[P(x_{13}),X_0]= (-6 x_{13} x_{12})\frac{\d }{\d x_{12}}+ (3 x_{13}^2)\frac{\d }{\d x_{13}}
+(-3 x_{13} x_{21}+x_{23} y_{13}+4 x_{23} y_{23})\frac{\d }{\d x_{21}}
 +(3 x_{13} x_{23})\frac{\d }{\d x_{23}}
+ (-3 x_{12} x_{21}+3 x_{13} x_{31}+3 x_{23} x_{32}+y_{13} y_{23}-2 y_{23}^2)\frac{\d }{\d x_{31}}
 +(2 x_{12} y_{13}+2 x_{12} y_{23})\frac{\d }{\d x_{32}}+ (6 x_{12} x_{23}+2 x_{13} y_{13}+2 x_{13} y_{23})\frac{\d }{\d y_{13}}
 +(3 x_{12} x_{23}+4 x_{13} y_{13}-5 x_{13} y_{23})\frac{\d }{\d y_{23}};
\end{dmath*}
\begin{dmath*}
X_3 =-[P(x_{21}),X_0]= (3 x_{12} x_{21}+3 x_{13} x_{31}-3 x_{23} x_{32}-y_{13}^2-y_{13} y_{23})\frac{\d }{\d x_{12}}
 + (-4 x_{23} y_{13}+2 x_{23} y_{23})\frac{\d }{\d x_{13}} + (3 x_{21}^2)\frac{\d }{\d x_{21}}
+ (-6 x_{23} x_{21})\frac{\d }{\d x_{23}}+ (3 x_{31} x_{21})\frac{\d }{\d x_{31}}
 + (-3 x_{21} x_{32}-5 x_{31} y_{13}+x_{31} y_{23})\frac{\d }{\d x_{32}}+ (-x_{21} y_{13}-4 x_{21} y_{23}-3 x_{23} x_{31})\frac{\d }{\d y_{13}}
 + (-5 x_{21} y_{13}-2 x_{21} y_{23}+3 x_{23} x_{31})\frac{\d }{\d y_{23}};
\end{dmath*}
\begin{dmath*}
X_4 =-[P(x_{23}),X_0]= (3 x_{12} x_{23}-4 x_{13} y_{13}-x_{13} y_{23})\frac{\d }{\d x_{12}}+ (-3 x_{13} x_{23})\frac{\d }{\d x_{13}}
 + (6 x_{23} x_{21})\frac{\d }{\d x_{21}}+ (-3 x_{23}^2)\frac{\d }{\d x_{23}}
 + (-2 x_{21} y_{13}-2 x_{21} y_{23})\frac{\d }{\d x_{31}}+ (3 x_{12} x_{21}-3 x_{13} x_{31}-3 x_{23} x_{32}+2 y_{13}^2-y_{13} y_{23})\frac{\d }{\d x_{32}}
 + (-3 x_{13} x_{21}+5 x_{23} y_{13}-4 x_{23} y_{23})\frac{\d }{\d y_{13}}+ (-6 x_{13} x_{21}-2 x_{23} y_{13}-2 x_{23} y_{23})\frac{\d }{\d y_{23}};
\end{dmath*}
\begin{dmath*}
X_5 =-[P(x_{31}),X_0]= (4 x_{32} y_{13}-2 x_{32} y_{23})\frac{\d }{\d x_{12}}
 + (-3 x_{12} x_{21}-3 x_{13} x_{31}+3 x_{23} x_{32}+y_{13}^2-3 y_{13} y_{23}+2 y_{23}^2)\frac{\d }{\d x_{13}}
+ (-3 x_{31} x_{21})\frac{\d }{\d x_{21}}
 + (5 x_{21} y_{13}-4 x_{21} y_{23}+3 x_{23} x_{31})\frac{\d }{\d x_{23}}
+ (-3 x_{31}^2)\frac{\d }{\d x_{31}}+ (6 x_{31} x_{32})\frac{\d }{\d x_{32}}
 + (6 x_{21} x_{32}-4 x_{31} y_{13}+2 x_{31} y_{23})\frac{\d }{\d y_{13}}
+ (3 x_{21} x_{32}-5 x_{31} y_{13}+7 x_{31} y_{23})\frac{\d }{\d y_{23}};
\end{dmath*}
\begin{dmath*}
X_6 =-[P(x_{32}),X_0]= (3 x_{32} x_{12})\frac{\d }{\d x_{12}}+ (4 x_{12} y_{13}-5 x_{12} y_{23}-3 x_{13} x_{32})\frac{\d }{\d x_{13}}
 + (2 x_{31} y_{13}-4 x_{31} y_{23})\frac{\d }{\d x_{21}}+ (3 x_{12} x_{21}-3 x_{13} x_{31}+3 x_{23} x_{32}-2 y_{13}^2+3 y_{13} y_{23}-y_{23}^2)\frac{\d }{\d x_{23}}
 + (-6 x_{31} x_{32})\frac{\d }{\d x_{31}}+ (3 x_{32}^2)\frac{\d }{\d x_{32}}+ (-3 x_{12} x_{31}-7 x_{32} y_{13}+5 x_{32} y_{23})\frac{\d }{\d y_{13}}
 + (-6 x_{12} x_{31}-2 x_{32} y_{13}+4 x_{32} y_{23})\frac{\d }{\d y_{23}});
\end{dmath*}
\begin{dmath*}
X_7 =-[P(x_{12}),X_6]=(-9 x_{12}^2)\frac{\d }{\d x_{13}}+ (-9 x_{12} x_{31}+3 x_{32} y_{13}+3 x_{32} y_{23})\frac{\d }{\d x_{21}}
 + (6 x_{12} y_{13}-3 x_{12} y_{23}+9 x_{13} x_{32})\frac{\d }{\d x_{23}}
+ (9 x_{32}^2)\frac{\d }{\d x_{31}}+ (18 x_{32} x_{12})\frac{\d }{\d y_{13}}+ (9 x_{32} x_{12})\frac{\d }{\d y_{23}};
\end{dmath*}
\begin{dmath*}
X_8 =-[P(x_{13}),X_4]=(9 x_{13}^2)\frac{\d }{\d x_{12}}+ (-9 x_{23}^2)\frac{\d }{\d x_{21}}
+ (9 x_{13} x_{21}-3 x_{23} y_{13}+6 x_{23} y_{23})\frac{\d }{\d x_{31}}
 + (-9 x_{12} x_{23}-6 x_{13} y_{13}+3 x_{13} y_{23})\frac{\d }{\d x_{32}}
+ (-9 x_{13} x_{23})\frac{\d }{\d y_{13}}+ (9 x_{13} x_{23})\frac{\d }{\d x_{23}};
\end{dmath*}
\begin{dmath*}
X_9 =-[P(x_{31}),X_3]=(9 x_{21} x_{32}-3 x_{31} y_{13}-3 x_{31} y_{23})\frac{\d }{\d x_{12}}+ (3 x_{21} y_{13}-6 x_{21} y_{23}-9 x_{23} x_{31})\frac{\d }{\d x_{13}}
 + (9 x_{21}^2)\frac{\d }{\d x_{23}}+ (-9 x_{31}^2)\frac{\d }{\d x_{32}}
+ (-9 x_{31} x_{21})\frac{\d }{\d y_{13}}+ (-18 x_{31} x_{21})\frac{\d }{\d y_{23}}.
\end{dmath*}
\end{dgroup*}
The hamiltonian $Q_b$ (satisfying (\ref{id1})):
\begin{dmath*}
Q_b= 6( (-4  b_2  b_4-12  b_5  b_7)  x_{21}  x_{13}+(12  b_0  b_7-8  b_2  b_3)  x_{23}  x_{13}+(-   b_0^2+4  b_1  b_4-4  b_2  b_5+4  b_3  b_6)  x_{13}  x_{31}+
(-4  b_2  b_6-12  b_3  b_8)  x_{13}  x_{32}+(-4  b_0  b_2-4  b_1  b_3+12  b_6  b_7)  y_{13}  x_{13}+(12  b_5  b_8+4  b_6^2)  x_{32}^2+
(4  b_0  b_6+8  b_1  b_5)  y_{23}  x_{32}+(   b_0^2-4  b_3  b_6)  y_{13}^2+(-   b_0^2-4  b_1  b_4+4  b_2  b_5+4  b_3  b_6)  y_{23}  y_{13}+
(   b_0^2-4  b_2  b_5)  y_{23}^2+(-12  b_1  b_9-4  b_4  b_6)  y_{23}  x_{31}+(4  b_5^2+12  b_6  b_9)  x_{31}^2+(4  b_3^2+12  b_4  b_7)  x_{23}^2+
(-4  b_0  b_3-4  b_2  b_4+12  b_5  b_7)  y_{23}  x_{23}+(-4  b_0  b_4+12  b_2  b_9-4  b_3  b_5)  y_{23}  x_{21}+(12  b_3  b_9+4  b_4^2)  x_{21}^2+
(-4  b_1  b_3-12  b_6  b_7)  y_{23}  x_{13}+(12  b_1  b_7+4  b_2^2)  x_{13}^2+(4  b_0  b_1+8  b_2  b_6)  y_{23}  x_{12}+(4  b_1^2+12  b_2  b_8)  x_{12}^2+
(-4  b_1  b_5-12  b_4  b_8)  x_{32}  y_{13}+(4  b_5  b_6-36  b_8  b_9)  x_{32}  x_{31}+(4  b_0  b_5+8  b_4  b_6)  y_{13}  x_{31}+
(-12  b_2  b_9-4  b_3  b_5)  x_{23}  x_{31}+(-   b_0^2+4  b_1  b_4+4  b_2  b_5-4  b_3  b_6)  x_{23}  x_{32}+(-4  b_2  b_4-12  b_5  b_7)  x_{23}  y_{13}+
(-4  b_1  b_3-12  b_6  b_7)  x_{12}  x_{23}+(-4  b_1  b_5-12  b_4  b_8)  x_{31}  x_{12}+(12  b_0  b_8-8  b_1  b_6)  x_{32}  x_{12}+
(-4  b_0  b_1-4  b_2  b_6+12  b_3  b_8)  x_{12}  y_{13}+(4  b_1  b_2-36  b_7  b_8)  x_{12}  x_{13}+
(-   b_0^2-4  b_1  b_4+4  b_2  b_5+4  b_3  b_6)  x_{12}  x_{21}+(4  b_3  b_4-36  b_7  b_9)  x_{21}  x_{23}+(12  b_0  b_9-8  b_4  b_5)  x_{31}  x_{21}+
(-12  b_1  b_9-4  b_4  b_6)  x_{32}  x_{21}+(4  b_0  b_4+8  b_3  b_5)  x_{21}  y_{13} ).
\end{dmath*}
The polynomial  $k_b=X_bC_2$ (see Theorem  \ref{cas}):
\begin{dmath*}
k_b=3 (b_0 x_{12} x_{21} y_{13}-b_0 x_{12} x_{21} y_{23}+3 b_0 x_{12} x_{23} x_{31}-3 b_0 x_{13} x_{21} x_{32}
-b_0 x_{13} x_{31} y_{13}+b_0 x_{23} x_{32} y_{23}-b_0 y_{13}^2 y_{23}
+b_0 y_{13} y_{23}^2-2 b_1 x_{12}^2 x_{21}
+4 b_1 x_{12} x_{13} x_{31}-2 b_1 x_{12} x_{23} x_{32}+2 b_1 x_{12} y_{13} y_{23}-2 b_1 x_{13} x_{32} y_{13}
+4 b_1 x_{13} x_{32} y_{23}
 -4 b_2 x_{12} x_{13} x_{21}+2 b_2 x_{12} x_{23} y_{13}+2 b_2 x_{12} x_{23} y_{23}
+2 b_2 x_{13}^2 x_{31}+2 b_2 x_{13} x_{23} x_{32}+2 b_2 x_{13} y_{13} y_{23}-2 b_2 x_{13} y_{23}^2
 +4 b_3 x_{12} x_{21} x_{23}-2 b_3 x_{13} x_{21} y_{13}-2 b_3 x_{13} x_{21} y_{23}-2 b_3 x_{13} x_{23} x_{31}-2 b_3 x_{23}^2 x_{32}+2 b_3 x_{23} y_{13}^2-2 b_3 x_{23} y_{13} y_{23}
 +2 b_4 x_{12} x_{21}^2
+2 b_4 x_{13} x_{21} x_{31}-4 b_4 x_{21} x_{23} x_{32}-2 b_4 x_{21} y_{13} y_{23}-4 b_4 x_{23} x_{31} y_{13}+2 b_4 x_{23} x_{31} y_{23}-2 b_5 x_{12} x_{21} x_{31}
-2 b_5 x_{13} x_{31}^2+4 b_5 x_{21} x_{32} y_{13}
-2 b_5 x_{21} x_{32} y_{23}+4 b_5 x_{23} x_{31} x_{32}-2 b_5 x_{31} y_{13} y_{23}+2 b_5 x_{31} y_{23}^2+2 b_6 x_{12} x_{21} x_{32}
 +2 b_6 x_{12} x_{31} y_{13}-4 b_6 x_{12} x_{31} y_{23}-4 b_6 x_{13} x_{31} x_{32}+2 b_6 x_{23} x_{32}^2-2 b_6 x_{32} y_{13}^2+2 b_6 x_{32} y_{13} y_{23}-6 b_7 x_{12} x_{23}^2
 +6 b_7 x_{13}^2 x_{21}-6 b_7 x_{13} x_{23} y_{13}+6 b_7 x_{13} x_{23} y_{23}-6 b_8 x_{12}^2 x_{31}+6 b_8 x_{12} x_{32} y_{13}
+6 b_8 x_{13} x_{32}^2+6 b_9 x_{21}^2 x_{32}
 -6 b_9 x_{21} x_{31} y_{23}-6 b_9 x_{23} x_{31}^2).
 \end{dmath*}
The polynomial  $M_b:=m_b-X_bQ_b$ (see Theorem  \ref{cas}):
\begin{dmath*}
M_b=6((-16 b_3^3-72 b_3 b_4 b_7+216 b_7^2 b_9) x_{23}^3+(-16 b_5^3-72 b_5 b_6 b_9
+216 b_8 b_9^2) x_{31}^3+(72 b_5 b_6 b_8+16 b_6^3-216 b_8^2 b_9) x_{32}^3
+(4 b_0 b_1 b_4-8 b_0 b_2 b_5+4 b_0 b_3 b_6+12 b_1 b_2 b_9-12 b_1 b_3 b_5+12 b_2 b_4 b_6-24 b_3 b_4 b_8
+12 b_5 b_6 b_7) y_{13}^3+(-12 b_0 b_1 b_4+12 b_0 b_3 b_6+36 b_1 b_2 b_9-12 b_1 b_3 b_5+12 b_2 b_4 b_6
-36 b_5 b_6 b_7) y_{23}^3+(72 b_3 b_4 b_9+16 b_4^3-216 b_7 b_9^2) x_{21}^3+(72 b_1 b_2 b_7+16 b_2^3
-216 b_7^2 b_8) x_{13}^3+(-16 b_1^3-72 b_1 b_2 b_8+216 b_7 b_8^2) x_{12}^3+(-48 b_0 b_5 b_6-48 b_1 b_5^2
-144 b_4 b_5 b_8-96 b_4 b_6^2) x_{31} x_{32} y_{13}+(48 b_0 b_5 b_6+96 b_1 b_5^2+144 b_1 b_6 b_9
+48 b_4 b_6^2) x_{31} x_{32} y_{23}+(-10 b_0^2 b_5-24 b_0 b_1 b_9-48 b_0 b_4 b_6-32 b_1 b_4 b_5
+40 b_2 b_5^2-24 b_2 b_6 b_9-56 b_3 b_5 b_6+72 b_3 b_8 b_9+24 b_4^2 b_8) x_{31} y_{13} y_{23}
+(10 b_0^2 b_6+48 b_0 b_1 b_5+24 b_0 b_4 b_8-24 b_1^2 b_9+32 b_1 b_4 b_6+56 b_2 b_5 b_6-72 b_2 b_8 b_9
+24 b_3 b_5 b_8-40 b_3 b_6^2) x_{32} y_{13} y_{23}+
(-12 b_0^2 b_3+48 b_2 b_3 b_5+48 b_3^2 b_6
+144 b_4 b_6 b_7) x_{12} x_{21} x_{23}+(6 b_0^2 b_5-72 b_0 b_1 b_9-24 b_2 b_5^2-216 b_2 b_6 b_9
-24 b_3 b_5 b_6-216 b_3 b_8 b_9-72 b_4^2 b_8) x_{12} x_{21} x_{31}+
\end{dmath*}
\begin{dmath*}
(-6 b_0^2 b_6+72 b_0 b_4 b_8
+72 b_1^2 b_9+24 b_2 b_5 b_6+216 b_2 b_8 b_9+216 b_3 b_5 b_8+24 b_3 b_6^2) x_{12} x_{21} x_{32}
+(-3 b_0^3-54 b_0 b_1 b_4+12 b_0 b_2 b_5+30 b_0 b_3 b_6+90 b_1 b_2 b_9-30 b_1 b_3 b_5-42 b_2 b_4 b_6
+36 b_3 b_4 b_8-162 b_5 b_6 b_7-324 b_7 b_8 b_9) x_{12} x_{21} y_{13}+(3 b_0^3+18 b_0 b_1 b_4
-12 b_0 b_2 b_5+6 b_0 b_3 b_6+18 b_1 b_2 b_9+42 b_1 b_3 b_5+30 b_2 b_4 b_6-36 b_3 b_4 b_8+54 b_5 b_6 b_7
+324 b_7 b_8 b_9) x_{12} x_{21} y_{23}+
(-9 b_0^3-18 b_0 b_1 b_4+36 b_0 b_2 b_5+90 b_0 b_3 b_6
+126 b_1 b_2 b_9-42 b_1 b_3 b_5+162 b_2 b_4 b_6-36 b_3 b_4 b_8-198 b_5 b_6 b_7
+324 b_7 b_8 b_9) x_{12} x_{23} x_{31}+(6 b_0^2 b_1-72 b_0 b_3 b_8-24 b_1^2 b_4-24 b_1 b_2 b_5
-216 b_2 b_4 b_8-216 b_5 b_7 b_8-72 b_6^2 b_7) x_{12} x_{23} x_{32}+(-6 b_0^2 b_2-48 b_0 b_1 b_3
+72 b_0 b_6 b_7-24 b_1 b_2 b_4+72 b_1 b_5 b_7+24 b_2^2 b_5-96 b_2 b_3 b_6+72 b_3^2 b_8
+216 b_4 b_7 b_8) x_{12} x_{23} y_{13}+
(-6 b_0^2 b_2-72 b_0 b_6 b_7-24 b_1 b_2 b_4-72 b_1 b_5 b_7
+24 b_2^2 b_5+48 b_2 b_3 b_6-72 b_3^2 b_8-216 b_4 b_7 b_8) x_{12} x_{23} y_{23}+(144 b_0 b_5 b_8
-48 b_1 b_5 b_6+144 b_4 b_6 b_8) x_{12} x_{31} x_{32}+
(-6 b_0^2 b_6-72 b_0 b_4 b_8-72 b_1^2 b_9
+48 b_1 b_4 b_6-24 b_2 b_5 b_6-216 b_2 b_8 b_9-72 b_3 b_5 b_8+24 b_3 b_6^2) x_{12} x_{31} y_{13}
+(12 b_0^2 b_6+48 b_0 b_1 b_5+48 b_1 b_4 b_6+48 b_2 b_5 b_6-48 b_3 b_6^2) x_{12} x_{31} y_{23}
+(54 b_0^2 b_8-24 b_0 b_1 b_6+24 b_1^2 b_5+72 b_2 b_5 b_8-24 b_2 b_6^2-144 b_3 b_6 b_8) x_{12} x_{32} y_{13}
+(-48 b_1^2 b_5-144 b_1 b_4 b_8+48 b_2 b_6^2+144 b_3 b_6 b_8) x_{12} x_{32} y_{23}+(26 b_0^2 b_1
+48 b_0 b_2 b_6-24 b_0 b_3 b_8-8 b_1^2 b_4-8 b_1 b_2 b_5-32 b_1 b_3 b_6-24 b_2 b_4 b_8+72 b_5 b_7 b_8
-120 b_6^2 b_7) x_{12} y_{13} y_{23}+(-144 b_0 b_4 b_7+48 b_2 b_3 b_4-144 b_3 b_5 b_7) x_{13} x_{21} x_{23}
+(-6 b_0^2 b_4+72 b_0 b_2 b_9+216 b_1 b_3 b_9+24 b_1 b_4^2+24 b_3 b_4 b_6+72 b_5^2 b_7
+216 b_6 b_7 b_9) x_{13} x_{21} x_{31}+(9 b_0^3-90 b_0 b_1 b_4-36 b_0 b_2 b_5+18 b_0 b_3 b_6
+198 b_1 b_2 b_9-162 b_1 b_3 b_5+42 b_2 b_4 b_6+36 b_3 b_4 b_8-126 b_5 b_6 b_7
-324 b_7 b_8 b_9) x_{13} x_{21} x_{32}+(6 b_0^2 b_3+72 b_0 b_5 b_7+24 b_1 b_3 b_4+216 b_1 b_7 b_9+
72 b_2^2 b_9-48 b_2 b_3 b_5-24 b_3^2 b_6+72 b_4 b_6 b_7) x_{13} x_{21} y_{13}+
(6 b_0^2 b_3+48 b_0 b_2 b_4
-72 b_0 b_5 b_7+24 b_1 b_3 b_4-216 b_1 b_7 b_9-72 b_2^2 b_9+96 b_2 b_3 b_5-24 b_3^2 b_6
-72 b_4 b_6 b_7) x_{13} x_{21} y_{23}+(6 b_0^2 b_3-72 b_0 b_5 b_7-24 b_1 b_3 b_4-216 b_1 b_7 b_9
-72 b_2^2 b_9-24 b_3^2 b_6-216 b_4 b_6 b_7) x_{13} x_{23} x_{31}+(-6 b_0^2 b_2+72 b_0 b_6 b_7
+24 b_1 b_2 b_4+216 b_1 b_5 b_7+24 b_2^2 b_5+72 b_3^2 b_8+216 b_4 b_7 b_8) x_{13} x_{23} x_{32}
+(-54 b_0^2 b_7+24 b_0 b_2 b_3+24 b_1 b_3^2-72 b_1 b_4 b_7-24 b_2^2 b_4+144 b_3 b_6 b_7) x_{13} x_{23} y_{13}
+(54 b_0^2 b_7-24 b_0 b_2 b_3+24 b_1 b_3^2+72 b_1 b_4 b_7-24 b_2^2 b_4-144 b_2 b_5 b_7) x_{13} x_{23} y_{23}
+(12 b_0^2 b_6-48 b_1 b_4 b_6-144 b_3 b_5 b_8-48 b_3 b_6^2) x_{13} x_{31} x_{32}+
(3 b_0^3-30 b_0 b_1 b_4
+36 b_0 b_2 b_5+6 b_0 b_3 b_6+18 b_1 b_2 b_9+42 b_1 b_3 b_5+30 b_2 b_4 b_6+108 b_3 b_4 b_8-90 b_5 b_6 b_7
+324 b_7 b_8 b_9) x_{13} x_{31} y_{13}+(36 b_0 b_1 b_4-36 b_0 b_3 b_6-108 b_1 b_2 b_9-12 b_1 b_3 b_5
+12 b_2 b_4 b_6+108 b_5 b_6 b_7) x_{13} x_{31} y_{23}+(6 b_0^2 b_1+48 b_0 b_2 b_6-72 b_0 b_3 b_8
-24 b_1^2 b_4+24 b_1 b_2 b_5+96 b_1 b_3 b_6-72 b_2 b_4 b_8-216 b_5 b_7 b_8-
72 b_6^2 b_7) x_{13} x_{32} y_{13}
+
\end{dmath*}
\begin{dmath*}
(-12 b_0^2 b_1-48 b_0 b_2 b_6+48 b_1^2 b_4-48 b_1 b_2 b_5-48 b_1 b_3 b_6) x_{13} x_{32} y_{23}
+(10 b_0^2 b_2+24 b_0 b_6 b_7+8 b_1 b_2 b_4-120 b_1 b_5 b_7-40 b_2^2 b_5-64 b_2 b_3 b_6-24 b_3^2 b_8
-72 b_4 b_7 b_8) x_{13} y_{13} y_{23}+(144 b_0 b_3 b_9+144 b_2 b_4 b_9-48 b_3 b_4 b_5) x_{21} x_{23} x_{31}
+
(12 b_0^2 b_4-144 b_1 b_3 b_9-48 b_1 b_4^2-48 b_2 b_4 b_5) x_{21} x_{23} x_{32}+
(48 b_0 b_3 b_4
+48 b_2 b_4^2+96 b_3^2 b_5+144 b_4 b_5 b_7) x_{21} x_{23} y_{13}+
(144 b_2 b_3 b_9+48 b_2 b_4^2
-48 b_3^2 b_5-144 b_4 b_5 b_7) x_{21} x_{23} y_{23}+(-144 b_0 b_6 b_9-144 b_1 b_5 b_9
+48 b_4 b_5 b_6) x_{21} x_{31} x_{32}+(144 b_1 b_4 b_9-144 b_2 b_5 b_9-48 b_3 b_5^2
+48 b_4^2 b_6) x_{21} x_{31} y_{13}+(-54 b_0^2 b_9+24 b_0 b_4 b_5+144 b_2 b_5 b_9+24 b_3 b_5^2
-72 b_3 b_6 b_9-24 b_4^2 b_6) x_{21} x_{31} y_{23}+(-12 b_0^2 b_5-48 b_0 b_4 b_6-48 b_1 b_4 b_5
+48 b_2 b_5^2-48 b_3 b_5 b_6) x_{21} x_{32} y_{13}+(6 b_0^2 b_5+72 b_0 b_1 b_9-48 b_1 b_4 b_5
-24 b_2 b_5^2+72 b_2 b_6 b_9+24 b_3 b_5 b_6+216 b_3 b_8 b_9+72 b_4^2 b_8) x_{21} x_{32} y_{23}
+(-26 b_0^2 b_4+24 b_0 b_2 b_9-48 b_0 b_3 b_5+24 b_1 b_3 b_9+8 b_1 b_4^2+32 b_2 b_4 b_5+8 b_3 b_4 b_6
+120 b_5^2 b_7-72 b_6 b_7 b_9) x_{21} y_{13} y_{23}+(-12 b_0^2 b_5+48 b_1 b_4 b_5+48 b_2 b_5^2
+144 b_2 b_6 b_9) x_{23} x_{31} x_{32}+(12 b_0^2 b_4+48 b_0 b_3 b_5-48 b_1 b_4^2+48 b_2 b_4 b_5
+48 b_3 b_4 b_6) x_{23} x_{31} y_{13}+
(-6 b_0^2 b_4+72 b_0 b_2 b_9-48 b_0 b_3 b_5+72 b_1 b_3 b_9
+24 b_1 b_4^2-96 b_2 b_4 b_5-24 b_3 b_4 b_6+72 b_5^2 b_7+216 b_6 b_7 b_9) x_{23} x_{31} y_{23}
+(36 b_0 b_1 b_4-36 b_0 b_3 b_6-108 b_1 b_2 b_9-12 b_1 b_3 b_5+12 b_2 b_4 b_6
+108 b_5 b_6 b_7) x_{23} x_{32} y_{13}+(-3 b_0^3-6 b_0 b_1 b_4+12 b_0 b_2 b_5-18 b_0 b_3 b_6
-54 b_1 b_2 b_9-30 b_1 b_3 b_5-42 b_2 b_4 b_6+36 b_3 b_4 b_8-18 b_5 b_6 b_7
-324 b_7 b_8 b_9) x_{23} x_{32} y_{23}+(-10 b_0^2 b_3-24 b_0 b_5 b_7-8 b_1 b_3 b_4+72 b_1 b_7 b_9
+24 b_2^2 b_9+64 b_2 b_3 b_5+40 b_3^2 b_6+120 b_4 b_6 b_7) x_{23} y_{13} y_{23}+(12 b_0^2 b_2
-144 b_1 b_5 b_7-48 b_2^2 b_5-48 b_2 b_3 b_6) x_{12} x_{13} x_{21}+(144 b_0 b_1 b_7-48 b_1 b_2 b_3
+144 b_2 b_6 b_7) x_{12} x_{13} x_{23}+(-12 b_0^2 b_1+48 b_1^2 b_4+48 b_1 b_3 b_6
+144 b_2 b_4 b_8) x_{12} x_{13} x_{31}+(-144 b_0 b_2 b_8+48 b_1 b_2 b_6-144 b_1 b_3 b_8) x_{12} x_{13} x_{32}
+(-48 b_1^2 b_3+144 b_1 b_6 b_7+48 b_2^2 b_6-144 b_2 b_3 b_8) x_{12} x_{13} y_{13}+
+(-48 b_0 b_1 b_2
-48 b_1^2 b_3-144 b_1 b_6 b_7-96 b_2^2 b_6) x_{12} x_{13} y_{23}+(-6 b_0^2 b_6+24 b_1 b_4 b_6+24 b_2 b_5 b_6
-216 b_2 b_8 b_9-72 b_3 b_5 b_8-24 b_3 b_6^2) x_{23} x_{32}^2+(10 b_0^2 b_3+24 b_0 b_2 b_4+24 b_0 b_5 b_7
+32 b_1 b_3 b_4-72 b_1 b_7 b_9-24 b_2^2 b_9-16 b_2 b_3 b_5-40 b_3^2 b_6-48 b_4 b_6 b_7) x_{23} y_{13}^2
+(-8 b_0^2 b_3-24 b_0 b_2 b_4+24 b_0 b_5 b_7+8 b_1 b_3 b_4-72 b_1 b_7 b_9+48 b_2^2 b_9-16 b_2 b_3 b_5
-16 b_3^2 b_6-48 b_4 b_6 b_7) x_{23} y_{23}^2+(-24 b_5^2 b_6-216 b_5 b_8 b_9-144 b_6^2 b_9) x_{31}^2 x_{32}
+(-24 b_0 b_5^2+72 b_1 b_5 b_9-48 b_4 b_5 b_6+216 b_4 b_8 b_9) x_{31}^2 y_{13}+(-72 b_0 b_6 b_9
-72 b_1 b_5 b_9+24 b_4 b_5 b_6) x_{31}^2 y_{23}+(144 b_5^2 b_8+24 b_5 b_6^2+216 b_6 b_8 b_9) x_{31} x_{32}^2
+(-8 b_0^2 b_5+24 b_0 b_1 b_9+32 b_1 b_4 b_5-16 b_2 b_5^2+24 b_2 b_6 b_9+32 b_3 b_5 b_6-72 b_3 b_8 b_9
+48 b_4^2 b_8) x_{31} y_{13}^2+
(10 b_0^2 b_5+24 b_0 b_1 b_9+24 b_0 b_4 b_6-16 b_1 b_4 b_5-40 b_2 b_5^2
-48 b_2 b_6 b_9+32 b_3 b_5 b_6-72 b_3 b_8 b_9-24 b_4^2 b_8) x_{31} y_{23}^2+(72 b_0 b_5 b_8-24 b_1 b_5 b_6
+72 b_4 b_6 b_8) x_{32}^2 y_{13}+
\end{dmath*}
\begin{dmath*}
(24 b_0 b_6^2+48 b_1 b_5 b_6-216 b_1 b_8 b_9-72 b_4 b_6 b_8) x_{32}^2 y_{23}
+(-10 b_0^2 b_6-24 b_0 b_1 b_5-24 b_0 b_4 b_8+24 b_1^2 b_9+16 b_1 b_4 b_6-32 b_2 b_5 b_6+72 b_2 b_8 b_9
+48 b_3 b_5 b_8+40 b_3 b_6^2) x_{32} y_{13}^2+
(8 b_0^2 b_6-24 b_0 b_4 b_8-48 b_1^2 b_9-32 b_1 b_4 b_6
-32 b_2 b_5 b_6+72 b_2 b_8 b_9-24 b_3 b_5 b_8+16 b_3 b_6^2) x_{32} y_{23}^2+(-5 b_0^3+2 b_0 b_1 b_4
+20 b_0 b_2 b_5+14 b_0 b_3 b_6-30 b_1 b_2 b_9+18 b_1 b_3 b_5-18 b_2 b_4 b_6-12 b_3 b_4 b_8-66 b_5 b_6 b_7
+108 b_7 b_8 b_9) y_{13}^2 y_{23}+(5 b_0^3+10 b_0 b_1 b_4-20 b_0 b_2 b_5-26 b_0 b_3 b_6-6 b_1 b_2 b_9
+18 b_1 b_3 b_5-18 b_2 b_4 b_6+12 b_3 b_4 b_8+102 b_5 b_6 b_7-108 b_7 b_8 b_9) y_{13} y_{23}^2
+(-24 b_1^2 b_2-216 b_1 b_7 b_8-144 b_2^2 b_8) x_{12}^2 x_{13}+(6 b_0^2 b_1+24 b_1^2 b_4-24 b_1 b_2 b_5
-24 b_1 b_3 b_6+72 b_2 b_4 b_8+216 b_5 b_7 b_8) x_{12}^2 x_{21}+(-216 b_0 b_7 b_8+24 b_1^2 b_3+72 b_1 b_6 b_7
+144 b_2 b_3 b_8) x_{12}^2 x_{23}+(18 b_0^2 b_8+24 b_1^2 b_5+72 b_2 b_5 b_8-72 b_3 b_6 b_8) x_{12}^2 x_{31}
+(-72 b_0 b_1 b_8+48 b_1^2 b_6+72 b_2 b_6 b_8+216 b_3 b_8^2) x_{12}^2 x_{32}+(24 b_0 b_1^2+72 b_0 b_2 b_8
+24 b_1 b_2 b_6-216 b_6 b_7 b_8) x_{12}^2 y_{13}+(-24 b_0 b_1^2-48 b_1 b_2 b_6+72 b_1 b_3 b_8
+216 b_6 b_7 b_8) x_{12}^2 y_{23}+(144 b_1^2 b_7+24 b_1 b_2^2+216 b_2 b_7 b_8) x_{12} x_{13}^2
+(-6 b_0^2 b_4-72 b_1 b_3 b_9-24 b_1 b_4^2+24 b_2 b_4 b_5+24 b_3 b_4 b_6-216 b_6 b_7 b_9) x_{12} x_{21}^2
+(18 b_0^2 b_7+24 b_1 b_3^2+72 b_1 b_4 b_7-72 b_2 b_5 b_7) x_{12} x_{23}^2+(-216 b_0 b_8 b_9+24 b_1 b_5^2
+144 b_1 b_6 b_9+72 b_4 b_5 b_8) x_{12} x_{31}^2+(72 b_0 b_6 b_8-72 b_1 b_5 b_8-48 b_1 b_6^2
-216 b_4 b_8^2) x_{12} x_{32}^2+(-8 b_0^2 b_1-24 b_0 b_2 b_6+24 b_0 b_3 b_8-16 b_1^2 b_4+8 b_1 b_2 b_5
-16 b_1 b_3 b_6-48 b_2 b_4 b_8-72 b_5 b_7 b_8+48 b_6^2 b_7) x_{12} y_{13}^2+
(-8 b_0^2 b_1+24 b_0 b_3 b_8
-16 b_1^2 b_4+32 b_1 b_2 b_5+32 b_1 b_3 b_6+24 b_2 b_4 b_8-72 b_5 b_7 b_8+48 b_6^2 b_7) x_{12} y_{23}^2
+(-18 b_0^2 b_7-72 b_1 b_4 b_7-24 b_2^2 b_4+72 b_3 b_6 b_7) x_{13}^2 x_{21}+
(72 b_0 b_2 b_7-72 b_1 b_3 b_7
-48 b_2^2 b_3-216 b_6 b_7^2) x_{13}^2 x_{23}+(-6 b_0^2 b_2+24 b_1 b_2 b_4-72 b_1 b_5 b_7-24 b_2^2 b_5
+24 b_2 b_3 b_6-216 b_4 b_7 b_8) x_{13}^2 x_{31}+(216 b_0 b_7 b_8-144 b_1 b_6 b_7-24 b_2^2 b_6
-72 b_2 b_3 b_8) x_{13}^2 x_{32}+(-72 b_0 b_1 b_7-24 b_0 b_2^2-24 b_1 b_2 b_3
+216 b_3 b_7 b_8) x_{13}^2 y_{13}+(72 b_0 b_1 b_7-24 b_1 b_2 b_3+72 b_2 b_6 b_7) x_{13}^2 y_{23}
+(216 b_0 b_7 b_9-144 b_2 b_3 b_9-24 b_2 b_4^2-72 b_4 b_5 b_7) x_{13} x_{21}^2+(-72 b_0 b_3 b_7+48 b_2 b_3^2
+72 b_2 b_4 b_7+216 b_5 b_7^2) x_{13} x_{23}^2+(6 b_0^2 b_5-24 b_1 b_4 b_5+24 b_2 b_5^2+72 b_2 b_6 b_9
-24 b_3 b_5 b_6+216 b_3 b_8 b_9) x_{13} x_{31}^2+(-18 b_0^2 b_8+72 b_1 b_4 b_8-72 b_2 b_5 b_8
-24 b_2 b_6^2) x_{13} x_{32}^2+(8 b_0^2 b_2+24 b_0 b_1 b_3-24 b_0 b_6 b_7-8 b_1 b_2 b_4+48 b_1 b_5 b_7
+16 b_2^2 b_5+16 b_2 b_3 b_6-48 b_3^2 b_8+72 b_4 b_7 b_8) x_{13} y_{13}^2+(-10 b_0^2 b_2-24 b_0 b_1 b_3
-24 b_0 b_6 b_7-32 b_1 b_2 b_4+48 b_1 b_5 b_7+40 b_2^2 b_5+16 b_2 b_3 b_6+24 b_3^2 b_8
+72 b_4 b_7 b_8) x_{13} y_{23}^2+(144 b_3^2 b_9+24 b_3 b_4^2+216 b_4 b_7 b_9) x_{21}^2 x_{23}
+(72 b_0 b_4 b_9-216 b_2 b_9^2-72 b_3 b_5 b_9-48 b_4^2 b_5) x_{21}^2 x_{31}+(-18 b_0^2 b_9+72 b_2 b_5 b_9
-72 b_3 b_6 b_9-24 b_4^2 b_6) x_{21}^2 x_{32}+(24 b_0 b_4^2-72 b_2 b_4 b_9+48 b_3 b_4 b_5
-216 b_5 b_7 b_9) x_{21}^2 y_{13}+
\end{dmath*}
\begin{dmath*}
(-72 b_0 b_3 b_9-24 b_0 b_4^2-24 b_3 b_4 b_5
+216 b_5 b_7 b_9) x_{21}^2 y_{23}+(-24 b_3^2 b_4-216 b_3 b_7 b_9-144 b_4^2 b_7) x_{21} x_{23}^2
+(-72 b_0 b_5 b_9+216 b_1 b_9^2+48 b_4 b_5^2+72 b_4 b_6 b_9) x_{21} x_{31}^2+
(216 b_0 b_8 b_9
-72 b_1 b_6 b_9-144 b_4 b_5 b_8-24 b_4 b_6^2) x_{21} x_{32}^2+(8 b_0^2 b_4-24 b_0 b_2 b_9
-24 b_1 b_3 b_9+16 b_1 b_4^2-32 b_2 b_4 b_5-32 b_3 b_4 b_6-48 b_5^2 b_7+72 b_6 b_7 b_9) x_{21} y_{13}^2
+
(8 b_0^2 b_4-24 b_0 b_2 b_9+24 b_0 b_3 b_5+48 b_1 b_3 b_9+16 b_1 b_4^2+16 b_2 b_4 b_5-8 b_3 b_4 b_6
-48 b_5^2 b_7+72 b_6 b_7 b_9) x_{21} y_{23}^2+(-216 b_0 b_7 b_9+72 b_2 b_3 b_9+24 b_3^2 b_5
+144 b_4 b_5 b_7) x_{23}^2 x_{31}+(6 b_0^2 b_3-24 b_1 b_3 b_4+216 b_1 b_7 b_9-24 b_2 b_3 b_5+24 b_3^2 b_6
+72 b_4 b_6 b_7) x_{23}^2 x_{32}+(-72 b_0 b_4 b_7+24 b_2 b_3 b_4-72 b_3 b_5 b_7) x_{23}^2 y_{13}
+(24 b_0 b_3^2+72 b_0 b_4 b_7+24 b_2 b_3 b_4-216 b_2 b_7 b_9) x_{23}^2 y_{23}+(18 b_0^2 b_9-72 b_1 b_4 b_9
+24 b_3 b_5^2+72 b_3 b_6 b_9) x_{23} x_{31}^2).
\end{dmath*}

\section{Acknowledgements}

The first author is indebted to Taras Skrypnyk for numerous stimulating discussions.

\end{document}